\documentclass[english]{article}
\usepackage{graphicx,psfrag,amsfonts}

\newcommand{\dist}{\mbox{$\,$ \rm dist \rm $\,$}}

\newtheorem{lemma}{Lemma}[section]
\newtheorem{theorem}[lemma]{Theorem}

\newtheorem{corollary}[lemma]{Corollary}
\newtheorem{Paragraph}[lemma]{}

\newtheorem{definition}[lemma]{Definition}
\newtheorem{remark}[lemma]{Remark}

\begin{document}

\title{Contractive piecewise continuous maps modeling networks of  inhibitory neurons.}
\author{ELEONORA CATSIGERAS \thanks{   Instituto de
Matem\'{a}tica y Estad\'{\i}stica Rafael Laguardia (IMERL),
 Fac. Ingenier\'{\i}a. Universidad de la Rep\'{u}blica. Uruguay.} \thanks{
  Corresponding author. Address: Herrera
  y Reissig 565. C.P. 11300. Montevideo. Uruguay. Tel-fax: +598-2711-0621. E-mail: eleonora@fing.edu.uy   }  ,
  \'{A}LVARO ROVELLA \thanks{Centro de Matem\'{a}tica. Fac. de Ciencias. Universidad
   de la Rep\'{u}blica. Uruguay.
  }  ,   RUBEN BUDELLI \thanks{Departamento de Biomatem\'{a}tica. Fac. de Ciencias. Universidad de
la Rep\'{u}blica. Uruguay}  }

\date{Revised version: December 17th., 2008.}

\maketitle

\begin{abstract}

 We prove that a topologically generic network (an open and dense set of networks)  of three or more
 inhibitory neurons have  periodic  behavior with a finite number of limit cycles that
persist under small perturbations of the structure of the network.

The network is modeled by the Poincar\'{e}  transformation which is piecewise continuous and locally contractive on a compact  region $B$ of  a finite
 dimensional  manifold, with the separation property: it transforms homeomorphically   the different
  continuity pieces of $B$
 into pairwise disjoint sets.

\end{abstract}

PACS 2008 codes: 87.19 lj, 87.19 ll, 87.19 lm

\vspace{.3cm}

MSC 2000 92B20, 34C25, 37G15, 24C28

\vspace{.3cm}


\section{Introduction}

We study the dynamics  of an abstract dynamical system modeling a
network composed with $n \geq 3$  neurons, where $n$ is arbitrarily large, that are reciprocally coupled by
inhibitory synapsis.

 Each neuron is modeled as a pacemaker of the type  integrate and fire
\cite{13}, including the case of relaxation oscillators.  The
internal variable $V_i = V_i(t)$ describing each neuron's potential, for $i = 1, 2 \ldots, n$ evolves
increasingly on time $t$ during the interspike intervals: $V_i(t)$ has positive first derivative $dV_i(t)/dt>0$ and negative second derivative $d^2V_i(t) / dt^2 <0 $ being the solution of a deterministic
autonomous differential equation of a wide general type.

When the potential $V_i$ reaches a given threshold
value,  the neuron $i$ produces an spike, its potential $V_i$ is reseted   to zero. It is suppose that all the neurons are inhibitory, there are no delays and the network is totally connected. When the neuron $i$  spikes, not only its potential $V_i$ changes, being reseted to zero, but also,
through the synaptical connections,  an action potential makes the
 other $n-1$ neurons $j \neq i$,
suddenly change their respective potentials $V_j$ with a jump of amplitude $-H_{ij} < 0$.

The instants of spiking are defined by the evolution of the system itself, and not predetermined by regular intervals of observation of the system. We analyze the state of the system immediately after each spike, in the sequence of  instants of spiking of the network. We prove that the state of the system after each spike is a function $F$ of the state after the prior spike. This function $F$ is the so called Poincar\'{e} map. The study of the dynamics by iteration of the Poincar\'{e} map  is not an artificial discretization of the real time dynamics. On the contrary, the dynamics of $F$ and its properties (for instance periodicity, chaotic attractors) are equivalent to those of the  system evolving in real time.

In \cite{MirolloSrogatz} the technique of the first return Poincar\'{e} map to a section transversal to the flux was first applied to study neuron networks, in that case, an homogeneous network of excitatory coupled pacemakers neurons.
In \cite{0}  the same technique is applied to networks of  inhibitory  cells,
analyzing the real time dynamics via a discrete Poincar\'{e} map $F$. This map is  locally contractive
and piecewise continuous in a compact set of ${\mathbb R}^{n-1}$. We include the proof of these properties
 in the section 2 of this paper.

In Section \ref{seccionperiodico}, we
prove  Theorem \ref{Teorema1}, which is the main abstract mathematical result:

\em Locally contractive piecewise continuous maps with the separation
property generically have   only periodic asymptotic behavior,
 with up to a finite number of
limit cycles that are persistent under small perturbations of the
map. \em
Generic  systems have a  topological meaning in this paper: they include an open and dense family of systems.

As a consequence we obtain the following applied result:

\em Generic neuron networks composed by $n \geq 3$ inhibitory
cells exhibit only periodic behavior  with  a finite number of  limit cycles
that are
 persistent under small perturbations of the set of parameter
values. \em

This is a  result generalizing the conclusions obtained for two
neurons networks in \cite{4} and \cite{6}.

On the other hand  non generic dynamics are structurable unstable: they are destroyed if the system is perturbed, even if the perturbation is arbitrarily small. We refer to those as bifurcating systems.

\vspace{.5cm}

 The  results of this paper are proved in an abstract and theoretical context, using the classical qualitative mathematical tools of the Topological Dynamical Systems Theory.
 The systems have discontinuities and evolve in finite but large dimension $n$. It is still mostly unknown, the dynamics  of  discontinuous and large dimensional  systems. That is why in sections 3 and 4 we include an abstract theory of topological dynamical systems with discontinuities. The piecewise continuity and the local contractiveness  help us  to obtain the thesis of persistent periodicity in the Lemma \ref{lema} of this paper. On the other hand, the separation property will play a fundamental role to obtain the thesis of density, and genericity  of the periodic behavior in the Theorem \ref{Teorema1}.

 Even being our mathematical analysis theoretically abstract, we observe that the conclusions about the dynamics of our model of neuron networks,  fit with those obtained by  experiments in computer simulations with mutually coupled identical neurons in  networks of up to $10^{10}$ cells, as reported in the following papers:

  In \cite{braun2}
 it was observed the transition among different periodic activity, indicating that  the simulation data perfectly fit to experimental and clinical observations. In the computer simulated experiments the alterations of the discharge patterns when passing from one periodic cycle to another, arise from changes of the network parameters, changes in the connectivity between cells, and also of external modulation.

  { In \cite{braun3} } the computer simulated experiment shows the dynamics of the network of a large number of coupled neurons. It was observed to be significantly different from the original dynamics of the individual cells: the system can be driven through different synchronization states.

  Our thesis of Theorem  \ref{Teorema1} is only applicable to deterministic systems. Nevertheless their conclusions also qualitatively fit
   with computer simulations of neural systems  with  randomness \cite{villa},  \cite{turovavilla}, which also  show the generation of  detectable preferred firing sequences.

\section{A mathematical model of the inhibitory neurons network.} \label{modelo}

We include the detailed proof of the mathematical translation from a physical
model of  $n$ inhibitory pacemaker neurons network to  the dynamics of iterations
of a piecewise
continuous contractive map $F:B \mapsto B$, locally contractive and with the separation
property, as first posed in \cite{0}. The model is applicable for any
finite number $n \geq 3$ of neurons in the network.

The phase space of the system is the compact cube $Q = [-1,1]^n \subset  \mathbb{R}^n$. A point in the phase space is $V = (V_1, V_2, \ldots, V_n)$, describing the potential $V_i$ of each of the neurons $i \in \{1,2, \ldots, n\}$. We assume that the the phase space is normalized: the threshold level of each of the neurons potentials is 1, so 1 is the maximum of $V_i$. Also the minimum $V_i$ is normalized to $-1$, and that the reset value of  $V_i  $, after a spike of the neuron $i$, is $0.$

\begin{definition}
\label{definicionModeloInicial}
\em {\bf The physical model. }
The point $V$ in the phase space $Q$ evolves on time $t$, during the interspike intervals of time, according to an autonomous differential equation and changes without delay in a discontinuous fashion in the exact spiking instants, according to a reseting-synaptical rule. The two regimes, during the interspike interval, and in the spiking instants respectively, are precisely defined according to the following assumptions:
\end{definition}

{\bf \ref{definicionModeloInicial}.1}  {\bf Inter-spike regime assumptions.}  $V_i(t)$ is the solution of a differential equation
\begin{equation} \label{ecuacionDiferencial}
\frac{dV_i}{dt} = \gamma _i (V_i), \; \; \gamma _i : [-1,1] \mapsto  \mathbb{R}, \; \;
\gamma _i \in {\cal C}^1, \; \; \gamma_i (V_i) >0, \; \gamma _i ' (V_i) < 0 \; \; \forall \, V_i \in [-1,1].
\end{equation}

where ${\cal C}^1$ denotes the space of  real functions in $[-1,1]$,  continuous and derivable with continuous  derivative in $[-1,1]$.

The assumption $\gamma_i >0$ reflects that each neuron potential in the inter-spike interval is strictly increasing while it does not receive interactions from the other neurons of the network. This comes from the hypothesis that each isolated neuron $i$ is of pacemaker type, i.e. from any initial state $V_i(0) \in [-1,1)$, the potential  spontaneously reaches the threshold level $1$ for some time $t = t_i >0$,  if none inhibitory synapsis is received in the time interval $[0, t_i]$.

The assumption $\gamma '_i < 0$, which we call the \em dissipative hypothesis \em reflects  that the cynetic energy $E_c = (1/2)(dV_i/dt)^2$  is decreasing on time while the potential freely evolves during the interspike intervals. In fact: $d E_c /dt = (dV_i/dt) (d^2V_i/dt^2) = \gamma_i (V_i) \gamma_i' (V_i) \gamma_i (V_i)  < 0. $

The most used example of this type of inter-spike evolution is the relaxation oscillator model of a pacemaker neuron, for which $\gamma_i(V_i) = - \alpha _i V_i  + \beta_i$ where $0 < \alpha _i < \beta _i$ are constants. For this type of cells the differential equation (\ref{ecuacionDiferencial}) is linear, and its solution can be explicitly written: $$V_i(t) = (\beta_i/ \alpha_i) -
 [(\beta_i/ \alpha _i)  - V_i(0)] exp (- \alpha _i t).
$$

{\bf \ref{definicionModeloInicial}.2}  {\bf Consequences of the inter-spike regime assumptions.} We define  the flux $$\Phi^t(V) = (\Phi_1^t(V_1), , \Phi_2 ^t(V_2)\ldots, \Phi_n^t(V_n))$$ as the solution with initial state $V = (V_1, V_2, \ldots, V_n)$ of the differential equations system given by (\ref{ecuacionDiferencial}). Precisely:
\begin{equation}
\label{ecuacionFlujo}
\frac{ d(\Phi_i^t(V_i))}{dt} = \gamma _i (\Phi_i^t(V_i))\; \; \; \forall t , \; \; \; \Phi_i^0(V_i) = V_i \end{equation} As $\gamma_i \in {\cal C}^1$ we can apply the general theory of differential equations to deduce the following results, as a consequence of the assumptions in (\ref{ecuacionDiferencial}):

\begin{itemize}
\item
Two different orbits by the flux do not intersect.

 \item If $B$ and $A$ are two $(n-1)$-dimensional topological and connected sub-manifolds of $\mathbb{R}^{n}$ transversal to the vector field $\gamma = (\gamma_1, \ldots, \gamma_n)$, then the
  flux transforms homeomorphically any set of initial states in $B$ onto its image set of final states in $A$.

  \item For each constant time $t$ it holds the Louville formula:
  \begin{equation}
  \label{formulaLiouville}
 \frac{ d(\Phi_i^t(V_i))}{dV_i} = exp \int _0^t \gamma' _i (\Phi_i^s(V_i)) \, ds.
  \end{equation}
\end{itemize}

{\bf \ref{definicionModeloInicial}.3}  {\bf Spiking instants computations.}
 For each initial state $V \in Q$ the first spiking instant $\overline t (V)$  in the network is defined as the first positive time such that at least one of the neurons of the network reaches the threshold level 1. This means that
    \begin{equation}
    \label{ecuacionTiempoSpike}
    \overline t (V) = \min _{1 \leq i \leq n} t_i(V_i), \; \;  \mbox{ where } \Phi_i^{t} (V_i) = 1 \; \Leftrightarrow \; \; t = t_i(V_i)
    \end{equation}
     \begin{equation}
        \label{formulaConjuntoJ}
        J(V) = \{ i \in \{1, 2, \ldots, n\} : \overline t (V) = t_i(V_i)\}\end{equation} is the set of neurons that reach the threshold level simultaneously at the instant $\overline t (V)$. It is standard to prove  that for an open and dense set of initial states there is a single neuron $i$ reaching the threshold level first, i.e. $\# J(V) = 1, \; \; J(V) = \{i\} $.

    \vspace{.2cm}

    {\bf \ref{definicionModeloInicial}.4}  {\bf Spiking-synaptical assumptions.} In the spiking instant $\overline t $ the reseting and inhibitory synaptical interaction without delay produces an instantaneous discontinuity $$\sigma : \Phi^{\overline t(V)} (V) \mapsto \sigma (\Phi^{\overline t(V)} (V)) $$ in the state of the system, according to the following formulae:

        \begin{itemize}
             \item If $\# J(V) = 1, \; \; \{i\} = J(V)$ then $\Phi_i^{\overline t(V)} (V_i) = 1$ and:
             \begin{eqnarray} \sigma ^i & = & (\sigma^i _1, \sigma^i_2, \ldots, \sigma^i _n) \nonumber \\
              \sigma^i_i (\Phi_i^{\overline t(V)} (V_i)) &= &0 \mbox{ (spiking-reseting rule) } \label{ecuacionSinapsis1} \\ \sigma^i_j (\Phi_j^{\overline t(V)} (V_j)) &= & \max \; \{-1, \; \Phi_j^{\overline t(V)} (V_j) - H_{ij} \} \; \; \forall \,  j \neq i \mbox{ (synaptic rule)} \label{ecuacionSinapsis2}  \end{eqnarray}

             where $H_{ij}>0$ is constant, depending only on $i, j$, and gives the instantaneous discontinuity jump in the potential of neuron $j \neq i$ produced through  the inhibitory synaptical connection from neuron $i$ to neuron $j$

             \item If $\# J(V) = k \geq 2, \; \; \{i_1, i_2, \ldots, i_k\} = J(V)$ then $\sigma$ is multiply defined, having $k$ possible vectorial values $\sigma ^{i_1}, \sigma^{i_2} , \ldots, \sigma ^{i_k}$, where $\sigma ^{i_h}$ is defined according to formulae (\ref{ecuacionSinapsis1})  and
(\ref{ecuacionSinapsis2}).

        \end{itemize}

        We also assume that the network whose nodes are the cells and whose sides are the synaptical inhibitory interactions $H_{ij}$, is a complete bidirectionally connected graph. Precisely:
       \begin{equation}
       \label{formulaEpsilo0GrafoCompleto}
       0 < \epsilon_0 = \min _{i \neq j} H_{ij}
       \end{equation}

 {\bf \ref{definicionModeloInicial}.5} {\bf Relative large dissipative assumption.  } We assume the following relations between the functional parameters $\gamma _i$ in the differential equations (\ref{ecuacionDiferencial}) governing the dissipative interspike regime, and the real parameters $H_{i,j}$ in the formula (\ref{ecuacionSinapsis2}) governing the spiking-synaptical regime.
    \begin{eqnarray}
    \max_{i \neq j} H_{ij} &< & \frac{1}{4} \label{ecuacionesLargeDissipativeness1} \\
    \max_{i, j} \; |\gamma_i (3/4) - \gamma _j (3/4)| &< & \frac{ \min _{i } \min _{  V_i \in [1/4,3/4]} |\gamma'_i(V_i)|}{4} \label{ecuacionesLargeDissipativeness2} \\
    \frac{\max_{i \neq j}{H_{ij}}}{\min_{i \neq j}{H_{ij}}} -1 & < & \frac{ \min _{i } \min _{  V_i \in [1/4,3/4]} |\gamma'_i(V_i)|}{4 \, \max_{ i} \; \gamma_i (3/4)} \label{ecuacionesLargeDissipativeness3}
    \end{eqnarray}

    Condition (\ref{ecuacionesLargeDissipativeness1}) assumes that the discontinuity synaptical jumps $H_{ij}$ are not relatively as large as the widest range $[0,1]$ of the potential of the cells when they act as oscillators between the reset value 0 and the threshold level 1, free of synpatical interactions.

    The hypothesis (\ref{ecuacionesLargeDissipativeness2}) and (\ref{ecuacionesLargeDissipativeness3}) verify for instance for homogeneous networks in which all the functions $\gamma_i$ and all the synaptic interactions $H_{ij}$ are constant independent of the neurons $i,j$. But as they are open conditions, they also verify if the network is not homogeneous but the neurons and the synaptical jumps are not very different. Finally they also verify for networks that are  very heterogeneous, but  the dissipative parameter of the system $ \min _{i } \min _{  V_i \in [1/4,3/4]} |\gamma'_i(V_i)|$ is large enough.

    The assumptions above can be also possed for some number $0 < a < 1/2$ instead of $1/4$ in inequality (\ref{ecuacionesLargeDissipativeness1}), the number  $1-a$ instead of  $3/4$ in the values of $V_i$, and $2/(1-2a)$ instead of the denominator $4$, of inequalities (\ref{ecuacionesLargeDissipativeness2}) and (\ref{ecuacionesLargeDissipativeness3}) . Nevertheless, and without loss of generality, in the computations of this work we will take the assumptions above with  $a = 1/4$ to fix the numerical bounds.

\begin{Paragraph} \label{comments}
{\bf Comments about the physical model.}
\end{Paragraph}

The hypothesis (\ref{ecuacionesLargeDissipativeness1}), (\ref{ecuacionesLargeDissipativeness2}) and (\ref{ecuacionesLargeDissipativeness3}) will allow to prove the so called separation property in Theorem \ref{teoremaSparationProperty} in this paper. This property is essential  to prove that the family of all the systems which exhibit a limit set formed only by a finite number of limit cycles is dense, which leads to the topological genericity of such systems.

We observe that the assumptions in (\ref{ecuacionDiferencial}), (\ref{ecuacionSinapsis1}) and (\ref{ecuacionSinapsis2}) are more general that what they a priori seem. In fact, if instead of the variables $V_i$ which describe the electric potentials of each of the neurons, we used other equivalent variables, the vector field $\gamma $ of the differential equation (\ref{ecuacionDiferencial}),  and the synaptical vectorial interaction $\sigma$ given by (\ref{ecuacionSinapsis1}) and (\ref{ecuacionSinapsis2}), would have other coordinate expressions.

For instance, each isolated cell $i$ acts as an oscilator, whose potential $V_i$ varies in the interval $[0,1]$. We can diffeomorphically change the variable $V_i$ to a new one $\widehat V_i \in [0,1]$,  called the \em phase \em of the oscilator, which by definition, evolves \em linearly \em with the time $t$, during a time constant $\tau_i$. In the new variables the differential equation governing the phase state $\widehat V_i$ will be $ d \widehat V_i / dt  = 1/\tau _i$ and the flux will be linear in $Q$.

In \cite{0} it is developed the model in such  phase variables $\widehat V_i$ for which the flux is linear, and it is defined the synaptical inhibitory interaction jumps $-s_{i,j} <0$ in the phase state $\widehat V_j $, when the phase $\widehat V_i$ reaches the threshold level 1. To be equivalent to the constant jumps $-H_{ij}$ in the old variables $V_j$, it is showed in \cite{0} that the  interaction jumps $-s_{i,j} <0$ in the new phase variables $\widehat V_j $, must be  functions $s_{i,j} (\widehat V_j)$, strictly increasing with $\widehat V_j$ and such that $\widehat V_j - s_{ij} (\widehat V_j)$ is also strictly increasing. In a widest model the functions $s_{i,j}(\widehat V_j)$ are continuous but not necessarily differentiable.

 In resume, up to a change of variables, the model assumed in this paper in hypothesis (\ref{ecuacionDiferencial}), (\ref{ecuacionSinapsis1}) and (\ref{ecuacionSinapsis2}), includes for instance the model in \cite {0} in which the  flux  is linear during the interspike interval regime, and the synaptic  jumps in the spiking instants adequately depend of the phase of the postsynaptic neuron.

\begin{definition}
\label{definicionModeloMatematico} {\bf The Mathematical Model. }
\em In this subsection we will define a Poincar\'{e} section $B \subset Q$ of the dynamical system modeling physically the network of $n$ inhibitory neurons defined in \ref{definicionModeloInicial}. We then shall define the first return  Poincar\'{e} map $F: B \mapsto B$. We will prove  that this map is piecewise continuous,  locally contractive and has the separation property. These properties justify the Definition \ref{definicionModeloAbstracto}, at the end of this section, in which we will model and analyze this kind of inhibitory neuron networks through  the abstract   mathematical discrete dynamical system defined by the iterates of its Poincar\'{e} map $F$.
\end{definition}

{\bf \ref{definicionModeloMatematico}.1. }{\bf The Poincar\'{e} section $B$.} Let $B \in Q = [-1,1]^n$ be the compact $(n-1)$-dimensional set defined as follows:
\begin{equation}
\label{ecuacionSeccionPoincare}
B = \bigcup_{k=1}^n \widehat B_k  \; \; \mbox{ where } \; \widehat B_k = \{V   \in Q: V_k = 0\}
\end{equation}
The topology in $B$ is defined in each $\widehat B_k$ as the induced  by its inclusion in the   $(n-1)$ dimensional subspace $\{V_k = 0\} $ of $\mathbb{R}^{n}$.  Each $\widehat B_k$ is transversal to the flux defined in \ref{definicionModeloInicial}.2 solution of the system of differential equations (\ref{ecuacionDiferencial}), because the vector field $\gamma$ in the second term of this differential equations has all its components strictly positive.

After each spike, the state of the system is in $B$, due to the reset rule in equality (\ref{ecuacionSinapsis1}). So the system returns infinitely many times to $B$ from any initial state $V \in Q$.

\vspace{.3cm}

\begin{figure}[h]
\psfrag{1}{$1 $}\psfrag{2}{  $2$}\psfrag{3}{$3$} \psfrag{a}{$a
$}\psfrag{b}{$b$} \psfrag{a}{$a $}\psfrag{b}{$b$} \psfrag{c}{$c
$}\psfrag{d}{$d$}\psfrag{e}{$e $}\psfrag{f}{$f$} \psfrag{g}{$g
$}\psfrag{h}{$h$}\psfrag{i}{$i $}\psfrag{j}{$j$} \psfrag{k}{$k
$}\psfrag{l}{$l$}\psfrag{m}{$m $}\psfrag{n}{$n$} \psfrag{o}{$o
$}\psfrag{p}{$p $}\psfrag{q}{$q$}
\begin{center}\caption{\footnotesize Model of a 3 neurons network in the 3-dimensional cube:
Reaching the threshold level of neurons 1,2 and 3 corresponds to
the front faces 1, 2 and 3 respectively of the cube. Points marked
in black correspond to the linear evolution from backward faces to
the front faces. Firing of neurons 1, 2 or 3
correspond to the jumping dotted lines from the respective front face to
its parallel backwards face. The figure shows the evolution
after 8 spikes of the neuronal system: 1,2,1,1,1,1,1 and 3.}
 \label{fig1}\includegraphics[scale=.5]{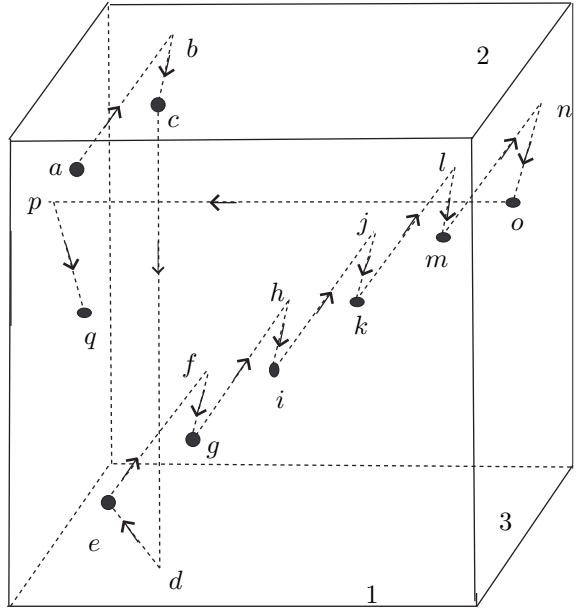}
\end{center}
\end{figure}

 {\bf The geometric illustration of the three neurons system. } In figure \ref{fig1} the cube $Q$ is represented for $n=3$ dimensions. The three coordinate axis $V_1, V_2, V_3$ are  the three edges of the cube that are hidden in dotted lines, at the rear part of the cube, intersecting pairwise orthogonally in the unique hidden vertix $O$ of the cube. The axis of $V_1 $ goes from $O$ to the front, the axis of $V_2$ to upwards, and the axis of $V_3$ to the right.

  The Poincar\'{e} section $B$ is formed by the three faces of the cube that correspond to at least one of the potentials $V_1, V_2, V_3$ equal to zero. $B$ is the union of the three faces of the cube at the rear part, that would not be seen if the cube were not transparent, each one  in a plane, orthogonal in the origin, to the  respective  coordinate axis.

 The increasing orbits of the flux inside the cube, are in that figure, parallel lines orthogonal to the plane of the figure. Due to the perspective each of the orbits is seen as a black dot in the figure. Each black dot, for instance ``a", represents a linear segment of an orbit. The dot ``a" is the orbit of the flux $\Phi^t(V)$ from the initial state $V= (V_1, V_2, 0) \in B$, in the left face of the cube (that would not be seen if the cube were not transparent),   to the front side of the cube, where the neuron  $i= 1$ reaches the threshold level 1.

 In that moment, the neuron $i=1$, whose potential increased to reach one (the state of the system is in the front vertical face of the cube), resets to zero, and the state of the system goes (through a dotted horizontal line parallel to the axis $V_1$), from the vertical front face of the cube (where $V_1 =1$) to the parallel vertical rear face in the Poincar\'{e} section $B$, where $V_1= 0$. (in the figure: from the point ``a" at front, to the point ``b" at back).

 In that  vertical rear face belonging to $B$, the other two neurons $j= 2, 3$, whose potentials $V_2, V_3$ were not reset, suffer a reduction of their potentials, of amplitudes $H_{12}$ and $H_{13}$ respectively, due to the inhibitory synaptic rule. That is why, the system does not stay in the point ``b" of the figure, but jumps to ``c", always in the rear face of the cube, corresponding to $V_1= 0$.

 At that instant, immediately after the first spike, from the point ``c"  in the backward rear face of the cube, the system starts  to evolve again according to the differential equation, in the inter-spike regime, moving on an orbit inside the cube.  In the figure this orbit corresponds to a segment orthogonal to the plane of the observer,   collapsed in the black dot ``c" due to the perspective. This new orbit arrives to the
 upper face of the cube, (also in the black dot ``c" of the figure), meaning that neuron $i=2$ arrived to the threshold level one.

  One could believe that  figure  \ref{fig1} is too particular, because the flux is linear inside the cube $Q$, with orbits that are parallel lines. But, due to the Tubular Flux Theorem, any flux tangent to a vector field $\gamma = (\gamma_1, \gamma_2, \gamma_3)$ such that $\gamma_i >0 $, after an adequate differentiable (generally non linear) change of variables $\xi $ in the space, becomes a linear flux whose orbits are parallel lines, and project orthogonally to a certain plane. One could transform all the dynamical system in these new coordinates $\xi (V)$, but in this case the synaptic jumps $H_{ij}$ in the old variables $V_j$ would become synaptic jumps $s_{ij}$ in the new variables $\xi_j (V)$. Maybe the matrix $(s_{ij})$ becomes dependent of the new variables $\widehat V_1= \xi_1(V),\; \; \widehat V_2 = \xi_2(V), \ldots, \widehat V_n = \xi_n(V)$, similarly to what was remarked in \ref{comments}.

  The Tubular Flux Theorem is also valid in any dimension $n \geq 2$, and this geometric model has all the data that we will develop analytically in this paper. In particular, the Tubular Flux Theorem and the linearization of the orbits, are used and analytically written in the proof of Theorem \ref{teoremaContraccion}.

The observer of the systems does not need to ``see" the orbits inside the 3-dimensional cube $Q$ of the figure \ref{fig1} to study its dynamics. He or she just see a point jumping inside the plane hexagon on which the cube $Q$ projects orthogonally to the flux, and orthogonally to the plane of the draw, from the viewpoint of the observer. There is a transformation $F$ from this hexagon to itself, giving the position of the point ``c" from the initial state ``a", then the position of ``e" from ``c", etc. The dynamics by iterates of this transformation $F$ describes exactly the same dynamics of the system, which indeed  evolves with continuous time $t$, but is disguised as discrete. It is not the system which was discrete nor the observer who made it discrete. The observer  positioned in the  adequate viewpoint, without modifying the system, just to see it in its   discrete disguise.

Nevertheless there is a problem with this discrete system, if two or more neurons got the threshold level simultaneously. For instance if the  point ``a" would be  in one of the three frontal sledges of the cube, that are sledges inside the hexagon, marked as not dotted lines in figure \ref{fig1}, there would be more than one possible consequent state $F(a)= ``c"$, depending on which frontal face the system chooses to  reset and apply the synaptic rule. That is why the transformation $F$ has discontinuities, and is not a uniquely defined map in  the frontal sledges, or lines of discontinuities.

All the arguments in this section could be obtained geometrically in the $n$-dimensional cube $Q$, just after its projection on an adequate poligon on an hyperplane, in which the state of the system evolves accordingly to a discrete transformation $F$.

\vspace{.2cm}

{\bf \ref{definicionModeloMatematico}.2. }{\bf The partition of $B$ in the continuity  pieces $B_i$. }  Recalling the definition of the spiking instant $\overline t(V)$ in equalities (\ref{ecuacionTiempoSpike}), and the definition of the set $J(V)$ of all the neurons that reach the threshold level at time $\overline t(V)$, in equality (\ref{formulaConjuntoJ}), we define the following subset $B_i$ of the Poincar\'{e} section $B$, for any $i = 1, 2 \ldots, n$:
\begin{equation}
\label{ecuacionBsubi}
B_i = \{V   \in B: \; i \in J(V)  \} = \{ V \in B: \; \overline t(V) = t_i(V_i)\}
\end{equation}
In other words, the set $B_i$ is formed by all the initial states $V$ in the Poincar\'{e} section $B$ such that the neuron $i$ reaches the threshold level before or at the same instant than all the other neurons of the network, from the initial state $V$.

From the implicit equation at right of formulae (\ref{ecuacionTiempoSpike}), we deduce that $B_i$  is compact, and that its interior $int(B_i)$ is formed by all the initial states for which $t_i(V_i) < t_j(V_j)$ for all $j \neq i$. Then $int(B_i) \bigcap int (B_j) = \emptyset \; \; \forall \, i \neq j$.

As the flux is strictly increasing inside $Q$, from any initial state $V \in B$ there exists a finite time $\overline t(V)$ defined by equalities (\ref{ecuacionTiempoSpike}). Therefore $V \in B_i$ for  some not necessarily unique $i\in \{1, 2, \ldots, n\}$. Then
$B = \bigcup _{i=1}^n B_i$ and the family of subsets
$ \{B_i\}_{i=1}^n$ is a topological finite partition of $B$ (i.e. it is a covering of $B$ with a finite number of compact sets whose interiors are pairwise disjoint.)

The compact sets $B_i$ are called \em continuity pieces. \em

We define the \em separation line $S$ \em of the partition $ \{B_i\}_{i=1}^n$, or \em line of discontinuities\em   as the union of the topological frontiers $\partial B_i$ of its subsets $B_i$. Precisely:
\begin{equation}
\label{formulaLineaSeparacionS} S = \bigcup_{i=1}^n \partial B_i = \bigcup_{i neq j} (B_i \cap B_j) = B \setminus \; \left ( \bigcup_{i=1}^n int B_i \right ) \end{equation}

{\bf \ref{definicionModeloMatematico}.3. }{\bf The first return Poincar\'{e} map $F$.}

The first return map $F: B \mapsto B$ to the Poincar\'{e} section $B =\bigcup_{i=1}^n B_i$  is  the finite collection of maps $f_i : B_i \mapsto B$  defined as $$f_i(V) = \sigma^i (\Phi^{\overline t(V)}(V)) \; \; \forall \, V \in B_i$$ where $\Phi$ is the solution flux defined in \ref{definicionModeloInicial}.2 of the system of differential equations (\ref{ecuacionDiferencial}), $\overline t(V)$ is the spiking instant defined by equalities (\ref{ecuacionTiempoSpike}) and $\sigma^i$ is the synaptic vectorial map defined in \ref{definicionModeloInicial}.4.

{\bf For simplicity we denote $F|_{B_i} = f_i$ and, when it is previously clear that $V \in B_i$, we  simply denote $F $ to refer to the uniquely well defined map $f_i$.}

We observe that $F$ is uniquely defined in $\bigcup_{i=1}^n int B_i$, and multi-defined in the separation line $S$.

Applying the formulae (\ref{ecuacionTiempoSpike}), (\ref{ecuacionSinapsis1}) and (\ref{ecuacionSinapsis2}), we deduce:
\begin{eqnarray}
F | _{B_i} (V) = f_i(V) &= &((f_i)_1, (f_i)_2, \ldots, (f_i)_n) \; \;  \forall \, V \in B_i  \mbox{ where } \nonumber \\
(F | _{B_i})_i(V) = (f_i)_i(V) &= & 0 \; \; = \; \; \max \; \{-1, \; \Phi_i^{ t_i(V_i)} (V_i) - H_{ii} \} \nonumber  \\
 (F | _{B_i})_j(V) = (f_i)_j(V) &= & \max \; \{-1, \; \Phi_j^{t_i(V_i)} (V_j) - H_{ij} \} \; \; \forall \,  j  \label{formulaReturnMapF}
\end{eqnarray}
where by convenience we agree to define $H_{ii} = +1$.

\vspace{.2cm}

{\bf \ref{definicionModeloMatematico}.4.} {\bf Piecewise continuity of the Poincar\'{e} map $F$.}

 The formula (\ref{formulaReturnMapF}) implies that $f_i = F | _{B_i}: B_i \mapsto B$ is continuous, and, as $B_i$ is compact, then $f_i(B_i)$ is also compact.

 The formula (\ref{formulaReturnMapF}) changes when one passes from  $B_i$ to $B_h$ with $i \neq h$, so $F$ is multidefined  in the points of $ S = \bigcup_{i \neq h} (B_i \cap B_h)$. Besides $F$ may be discontinuous in $V^0 \in B_i \cap B_h$ because $$\lim_{V \in int B_i, V \rightarrow V^0} F(V) = f_i(V^0)$$ is not necessarily equal to $$    \lim_{V \in int B_h, V \rightarrow V^0} F(V) = f_h(V^0) = $$.

  \vspace{.2cm}
  {\bf Remark \ref{definicionModeloMatematico}.4: } We agree to define the image set  $F(V)$ of a point $V \in B$ as $\{f_i(V): i \mbox{ such that } V \in B_i\}$. The image set $F(V)$ is a single point if $V \in int (B_i)$ because $int B_i$ does not intersect $B_j$ for $j \neq i$. The image set $F(A)$ of a set $A \subset B$ is by definition
  $F(A) = \bigcup _{V \in A} F(V)$.

  \vspace{.2cm}

  {\bf \ref{definicionModeloMatematico}.5.} {\bf The positive  reduced Poincar\'{e} section $B^+$.} We will define a subset $B^+ \subset B$ such that $F^p(B) \subset B^+$ for all $p$ large enough. Our aim is to study the limit set of the orbits, therefore the last property  allows us to restrict $F$ to $B^+$.
  \begin{equation} \label{ecuacionSeccionBmas}
  B^+ = \{V \in B: \; 0 \leq V_i \leq 1- \epsilon_0\; \; \forall \, i = 1, 2, \ldots, n  \}
  \end{equation}
  where $\epsilon _0 >0$ is the minimum of the absolute values of  the synaptic interactions $H_{ij}>0$ for $i \neq j$, as assumed in \ref{definicionModeloInicial}.4, equality (\ref{formulaEpsilo0GrafoCompleto}). Also, by hypothesis (\ref{ecuacionesLargeDissipativeness1}) we have \begin{equation}
    \label{epsilonCeroMenorUnCuarto} 0 < \epsilon_0 < \frac{1}{4}\; \; \; , \ \ \ \ \ \ \ \ \ \ \; \; \; \; \frac{3}{4} < 1 - \epsilon_0 < 1
  \end{equation}

  \vspace{.2cm}

  {\bf \ref{definicionModeloMatematico}.6. } {\bf Properties of the positive  reduced Poincar\'{e} section $B^+$.}

  The set $B^+$ is \em homeomorphic to a compact ball in $\mathbb{R}^{n-1}$, \em (property that the whole Poincar\'{e} section $B$ does not have). In fact, $B^+$ is the union of $n$ compact squares $\widehat B^+_k = \{V \in \mathbb{R}^n:\; V_k= 0 , \; 0 \leq V_i \leq 1 - \epsilon_0 \} = \{0\} \times [0, 1-\epsilon_0]^{n-1} \subset \mathbb{R}^{n-1}$ such that, for $h \neq k$: $\widehat B^+_k \cap \widehat B^+ _h \neq \emptyset$ is formed only by the $(n-2)$-dimensional lines $\{V_k= V_h= 0\}$ in the frontiers of both squares.

  \vspace{.5cm}

  The  dynamics properties of $F$ restricted to the positive Poincar\'{e} section $B^+$, which precisely justify the restriction to $B^+$, will be stated and proved in Theorems \ref{teoremaBmasinvariante}, \ref{teoremaFinjective}, \ref{teoremaSparationProperty} and \ref{teoremaContraccion}. They justify the definition of the abstract mathematical model in \ref{definicionModeloAbstracto}, whose dynamics in the future and attractors will be studied in the following sections.

  \vspace{.2cm}

  {\bf \ref{definicionModeloMatematico}.7. } {\bf The positive continuity pieces $B_i^+$ of the Poincar\'{e} map.}
  We define$$B_i^+ = B^+ \cap B_i$$ where $B^+$ is the positive reduced Poincar\'{e} section defined in \ref{definicionModeloMatematico}.5 and $B_i$ are the continuity pieces of the Poincar\'{e} map $F$, defined in \ref{definicionModeloMatematico}.2, equality (\ref{ecuacionBsubi}), and in \ref{definicionModeloMatematico}.4.

  \vspace{.2cm}

  {\bf Remark: } It is rather technical to prove that  \em $B_i^+$ is homeomorphic to a compact ball in ${\mathbb{R}^{n-1}}$.  \em We sketch here a proof, leaving the technical details:   $B_i^+  $ is the pre-image in $B^+$ by the flux $ \Phi = \Phi ^{\overline t(\cdot)} (\cdot)$, of the $(n-1)$ square $A_i = \{  V \in Q: V_i = 1   \}$. The flux is injective and continuous from $B^+$ onto its image $\Phi(B^+) \subset A= \bigcup_{A_i}$, because it is transversal to $A$ and to $B^+$ and two different orbits of the flux do not intersect. Any two  orbits that  intersect $B^+$ in different points do intersect  $A$ in different points. Continuous and injective maps $\Phi$ from a (homeomorphic) ball $B^+$ in ${\mathbb{R}^{n-1}}$ onto a set $\Phi(B^+) \subset \mathbb{R}^{n-1}$, are homeomorphisms, due to the Theorem of the Invariance of the Domain. Therefore $\Phi(B^+) $ is homeomorphic to a $(n-1)$-dimensional compact ball, contained in $A$.

   Also $\Phi(B^+) \cap A_i$ is. To prove this last assertion, be aware that the intersection of two compact homeomorphic balls is not necessarily a single homeomorphic ball, but it holds in our model, because the frontiers of $\Phi(B^+)$  and $A_i$ have some symmetric properties due to the fact that the flux $\Phi$ has $n$ components $\Phi_k$, each one depending only on the respective single variable $V_k$.

   Then, $B^+_i$ is the homeomorphic image by $\Phi^{-1}$ of the compact (homeomorphic) ball $\Phi(B^+) \cap A_i$. $\; \; \Box$
  \begin{theorem}.
  \label{teoremaBmasinvariante} {\bf The return map to the positive Poincar\'{e} section $B^+$. }

  The positive reduced Poincar\'{e} section $B^+ \subset B$ defined in \ref{definicionModeloMatematico}.5, is forward invariant by the Poincar\'{e} map $F|_{B^+}: B^+ \mapsto B^+$, and it is  reached  from any initial state in $B$. Even more,

   $F(B^+) \subset \; \bigcup_{i=1}^n \; \{ V \in B: V_i= 0, \; \; \; 0 < V_j \leq 1 - \epsilon_0\; \; \forall  \, j \neq i   \}\; \subset B^+,$

   and there exists $p \geq 1$ such that $F^p (B) \subset B^+.$

  \em (Recall that the constant $\epsilon_0 >0$ defined in Equality (\ref{formulaEpsilo0GrafoCompleto}) verifies the hypothesis (\ref{epsilonCeroMenorUnCuarto}).) \em

  \end{theorem}

  To prove Theorem \ref{teoremaBmasinvariante} we will use the following lemma:
  \begin{lemma}
  \label{lemaTiempoMinimo} There exists a constant positive minimum time $T $ $$T=
  \frac{\epsilon_0}{ \max_k\, \gamma_k(3/4)} >0$$ such that,
   if $V \in B$ verifies $ V_i \leq 1 - \epsilon_0 \; \; \forall \, 1 \leq i \leq n$, then the interspike interval $\overline t (V) \geq T$.
  \end{lemma}

  {\em Proof:} According to the formula (\ref{ecuacionTiempoSpike}): $ \; \; \overline t(V) = \min_{ i  } t_i(V_i)$ where $t= t_i(V_i) $ is the solution of the implicit equation $\Phi_i^{t}(V_i) = 1$.
  We integrate the differential equation (\ref{ecuacionDiferencial}) with initial condition $V_i$, and recall that $\gamma_i (V_i) >0$, while the real solution  $\Phi_i^s(V_i) \leq 1$  is strictly increasing with $s$ (for $V_i$ constant) and it is the solution of an autonomous differential equation. Using the hypothesis $ V_i \leq 1 -\epsilon_0$, and applying  the inequality (\ref{epsilonCeroMenorUnCuarto}), we obtain:
  \begin{equation}
   \label{EcuacionDiferencialIntegrada}
   \Phi_i^{ t } (V_i) = V_i + \int _0 ^t  \frac{d \Phi^s_i(V_i)}{ds} \; ds = V_i + \int_0^t \gamma_i(\Phi_i^s(V_i))\, ds \end{equation}
  $$1 = \Phi_i ^{t_i(V_i)} =  \Phi_i ^{t_i (V_i) - \underline {t_i}(V_i)} \left (\Phi^{\underline{t_i}(V_i)} (V_i)      \right )     = \Phi_i ^{t_i (V_i) - \underline {t_i}(V_i)} \left (1-\epsilon_0    \right )$$ $$1 = 1- \epsilon _0 + \int _{0 }^{t_i(V_i) - \underline {t_i} (V_i)}\gamma_i(\Phi_i^s(1-\epsilon_0))\, ds$$
  where $0 \leq \underline {t_i}(V_i) < t_i(V_i)$  and $\Phi^{\underline{t_i}(V_i)} (V_i)  = 1- \epsilon_0$, being $ \underline {t_i} (V_i)$  the  time that takes the flux $\Phi_i^t(V_i)$ to be equal to $1 - \epsilon_0$ from the initial state $V_i \leq 1 -\epsilon_0 $.

  Recall that $\gamma_i (V_i) >0$ is strictly decreasing with $V_i$:

  $$ \gamma_i(\Phi_i^s(1-\epsilon_0)) \leq \gamma_i (\Phi_i^0(1-\epsilon_0))= \gamma_i(1- \epsilon_0) < \gamma_i(3/4) \; \; \; \forall s \geq 0$$
  $$1 \leq  1- \epsilon _0 + \int _0^{t_i(V_i) - \underline {t_i} (V_i) } \gamma_i(3/4) \, ds $$

  $$\epsilon_0  \leq \gamma_i (3/4) \; \; [t_i(V_i) - \underline {t_i} (V_i)] \leq \gamma_i(3/4) \;  t_i(V_i) $$ $$ \Rightarrow \; \; t_i(V_i) \geq T = \frac{\epsilon_0}{ \max_k\, \gamma_k(3/4)} \; \; \; \forall \, i , \; \; \Rightarrow \; \; \overline t(V) = \max_{i} t_i(V_i) \geq T. \;\; \Box$$

  \vspace{.3cm}

  {\em \bf Proof of Theorem \ref{teoremaBmasinvariante}:}  It is enough to prove the following two assertions:

  \vspace{.3cm}

  {\bf Assertion \ref{teoremaBmasinvariante}.A:} $F_j(V) \leq 1 - \epsilon_0 \; \; \; \forall \, V \in B $ (even if $V \not \in B^+$),  $ \; \forall \, j = 1, 2 \ldots, n.$

  \vspace{.3cm}
  {\bf Assertion \ref{teoremaBmasinvariante}.B:} There exists a constant $\epsilon_1 >0$ such that for all $V \in B_i$, if $F_j(V) \leq 0$ for some  $j \neq i$, then $F_j(V) - V_j \geq \epsilon_1$.

  \vspace{.3cm}

  Note that the assertion \ref{teoremaBmasinvariante}.B states its thesis in particular if $V \not \in B^+$, and also if $V \in B^+$ and $V_j= 0$.

  Recall that from the formulae (\ref{formulaReturnMapF}) of the Poincar\'{e} map $F$: $F_i(V) = 0$ for all $V \in B_i$. Observe that, being $V_j \geq -1$ for all $V \in B$, from the Assertion  \ref{teoremaBmasinvariante}.B  we deduce that the first number $p \geq 1$ of iterates of $F$ such that $F^p(B) \subset B^+$ is at most equal to $1 + \mbox{Integer-Part} (1/\epsilon_1)$.
\vspace{.3cm}

To prove the Assertion \ref{teoremaBmasinvariante}.A, apply the formulae (\ref{formulaReturnMapF}) of the return Poincar\'{e} map $F$, and recall the assumptions (\ref{formulaEpsilo0GrafoCompleto}), (\ref{ecuacionesLargeDissipativeness1}). If $V \in B_i$ then
\begin{equation}
\label{formulaReturnMapF2}
F_i(V) = 0, \; \; \; F_j (V) = \max \{-1, \; \; \Phi_j^{ \overline t (V)}(V_j) - H_{ij} \} \leq \; 1- \min_{i \neq j} H_{ij} = 1- \epsilon_0\end{equation}

To prove the Assertion \ref{teoremaBmasinvariante}.B, fix $V \in B_i$ such that, for some $j \neq i$ \begin{equation}
    \label{equationefejotanopositivo}
F_j(V) \leq 0 \end{equation}

Use the formulae (\ref{formulaReturnMapF2}). We assert that
\begin{equation}
\label{ecuacionPhiMenorUnCuarto}
\Phi_j^{ \overline t (V)}(V_j) < \frac{1}{4}\end{equation}
In fact, if it were greater or larger than $1/4$, as $H_{ij} < 1/4$ due to hypothesis (\ref{ecuacionesLargeDissipativeness1}),  the formulae (\ref{formulaReturnMapF2}) would imply that $F_j(V_j) >0$ contradicting our hypothesis (\ref{equationefejotanopositivo}).

Due to the hypothesis of the differential equation (\ref{ecuacionDiferencial}), the function $\gamma_j(V_j) $ is strictly decreasing with $V_j$, and the flux $\Phi_j^t$ is strictly increasing with $t$. Use the integrate expression (\ref{EcuacionDiferencialIntegrada}) of the differential equation, to compute $\Phi^{\overline t(V)} _j (V_j)$, the inequality (\ref{ecuacionPhiMenorUnCuarto}) and the Lemma \ref{lemaTiempoMinimo}, to deduce:
$$0 \leq s \leq \overline t(V) \; \; \Rightarrow \; \; \; V_j \leq\Phi_j ^s(V_j) \leq \Phi_j^{\overline t(V)}(V_j) < \frac{1}{4}\; \; \; \Rightarrow \; \; \; \gamma_j(\Phi^s_j(V_j)) > \gamma_j(1/4)$$
$$\Rightarrow \; \; \; \int_0^{\overline t(V)} \gamma_j(\Phi_j^s(V_j))j \; ds > \gamma_j(1/4) \cdot \overline t(V) \geq \min_{k} \gamma _k(1/4) \cdot T = \frac{\min _{k} \gamma _k(1/4) \, \epsilon_0}{\max_{k} \gamma_k (3/4)}$$

Recalling the integral equation (\ref{EcuacionDiferencialIntegrada}) and the formula (\ref{formulaReturnMapF2}) of the return map $F$, we deduce:
$$F_j(V) - V_j \geq \Phi^{\overline t(V)} (V_j) - V_j - H_{ij} =   \int_0^{\overline t(V)} \gamma_j(\Phi_j^s(V_j))j \; ds - H_{ij} $$
$$F_j(V) - V_j \geq \epsilon_0 \left(\frac{\min _{k} \gamma _k(1/4) }{\max_{k} \gamma_k (3/4)} - \frac{\max_{i \neq j}{H_{ij}}}{\epsilon_0}    \right ) = \epsilon_1$$

To end the proof it is enough to show that $\epsilon_1 >0$. Recall from equality (\ref{formulaEpsilo0GrafoCompleto}) that $\epsilon_0 = \min_{i \neq j} H_{ij}>0$

$$\frac{\epsilon_1}{\epsilon_0}  =\frac{\min _{k} \gamma _k(1/4) }{\max_{k} \gamma_k (3/4)} - \frac{\max_{i \neq j}H_{ij}}{\min_{i \neq j}  H_{ij}} = \frac{\min _{k} \gamma _k(1/4)-\max_{k} \gamma_k (3/4)}{\max_{k} \gamma_k (3/4)} - \left (\frac{\max_{i \neq j}H_{ij}}{\min_{i \neq j}  H_{ij}} -1 \right )  $$
$$\frac{\epsilon_1}{\epsilon_0}= \frac{ \gamma _h(1/4)- \gamma_h(3/4) + \gamma_h(3/4) - \gamma_k (3/4)}{\max_{k} \gamma_k (3/4)} - \left (\frac{\max_{i \neq j}H_ {ij}}{\min_{i \neq j}  H_{ij}} -1 \right ) $$
where we have taken $h$ and $k$ such that $\gamma_h(1/4) = \min _k \gamma_k(1/4), \; \; \; \gamma_k (3/4) = \max_k \gamma_k(3/4)$.
$$\frac{\epsilon_1}{\epsilon_0}\geq  \frac{ \gamma _h(1/4)- \gamma_h(3/4) -| \gamma_h(3/4) - \gamma_k (3/4)|}{\max_{k} \gamma_k (3/4)} - \left (\frac{\max_{i \neq j}H_ {ij}}{\min_{i \neq j}  H_{ij}} -1 \right ) $$
$$\frac{\epsilon_1}{\epsilon_0}\geq  \frac{ \gamma _h(1/4)- \gamma_h(3/4) -\max_{h \neq k}| \gamma_h(3/4) - \gamma_k (3/4)|}{\max_{k} \gamma_k (3/4)} - \left (\frac{\max_{i \neq j}H_ {ij}}{\min_{i \neq j}  H_{ij}} -1 \right ) $$
Applying the mean value theorem of the derivative   of $\gamma_h$, which is negative due to the dissipation hypothesis of the differential equation in assumption (\ref{ecuacionDiferencial}), we obtain  $$0 < \gamma_h(1/4) - \gamma_h(3/4) = \left(\frac{3}{4} -\frac 1 4 \right) \cdot \left . (-{\gamma ' } _h(\chi)) \right |_{\chi \in [1/4, 3/4]}  \geq \frac{\min _i \, \min_{V_i \in [1/4, 3/4]} |\gamma'_i(V_i)|}{2}  $$
The last inequalities and the assumptions (\ref{ecuacionesLargeDissipativeness2}) and (\ref{ecuacionesLargeDissipativeness3}) of relative large dissipative, imply:
$$\frac{\epsilon_1}{\epsilon_0}>  \frac{\min _i \, \min_{V_i \in [1/4, 3/4]} |\gamma'_i(V_i)|}{ \max_{k} \gamma_k (3/4)} \left (\frac{1}{2} - \frac{1}{4} - \frac{1}{4}  \right ) = 0
\; \; \Box $$

  \vspace{.5cm}

  \begin{remark} \em
  {\bf Formula of the Poincar\'{e}  map in $B^+$.}

  As a consequence of  Theorem \ref{teoremaBmasinvariante}, from now on we will restrict the Poincar\'{e} map $F$ to the positive section $B^+$. In fact, from the statements of Theorem \ref{teoremaBmasinvariante} it is deduced that the forward dynamics and the limit set of the orbits to the future, of the restricted $F$, will be the same as those of $F$ in the whole Poincar\'{e} section.

  Due to Theorem \ref{teoremaBmasinvariante} if $V \in B^+$ then $F(V) \subset B^+$. Therefore
  $(F|_{B_i})_j (V) \geq 0 \; \; \forall \, i, j$. Using the formula (\ref{formulaReturnMapF}) we can rewrite the expression of the Poincar\'{e} map, without the maximum:
  \begin{eqnarray}
F | _{B_i} (V) = f_i(V) &= &((f_i)_1, (f_i)_2, \ldots, (f_i)_n) \; \;  \forall \, V \in B_i^+  \mbox{ where } \nonumber \\
(F | _{B_i^+})_i(V) = (f_i)_i(V) &= & 0 \; \; = \; \;  \; \Phi_i^{ t_i(V_i)} (V_i) - H_{ii}  \nonumber  \\
 (F | _{B_i^+})_j(V) = (f_i)_j(V) &= &   \Phi_j^{t_i(V_i)} (V_j) - H_{ij}  \; \; \forall \,  j  \label{formulaReturnMapF3}
\end{eqnarray}
  \end{remark}

  \begin{theorem}
  \label{teoremaFinjective} {\bf Local injectiveness of the Poincar\'{e} map.}

  The Poincar\'{e} map $F $ defined in formulae (\ref{formulaReturnMapF3}), restricted to each of its positive continuity pieces $B^+_i $ defined in \em \ref{definicionModeloMatematico}.7,  \em is injective.
  \end{theorem}
  {\em Proof: }  Fix a continuity piece $B_i^+$ of $F$ in the positive Poincar\'{e} section $B^+$. The piece $B_i^+$ will remain  fixed along this proof. Therefore we will denote $F$ instead of $f_i$.

  Take $V, W \in B_i^+$.such that $F(V) = F(W) \in B^+$. We must  prove that $V = W$.

   Due to the formulas (\ref{formulaReturnMapF3}) of the Poincar\'{e} map: $F_i(V) = F_i(W) = 0$
  and
  $$F_j(V) = \Phi_j^{\overline t(V)}(V_j) + H_{ij} = \Phi_j^{\overline t (W)} (W_j) + H_{ij} = F_j(W) \; \; \ \forall \, j $$

As $H_{ij}$ is constant, we deduce that
$$\Phi_j ^{\overline t(V)}(V_j) = \Phi_j^{\overline t (W)} (W_j) \; \; \; \forall \, j $$

Therefore the vectorial flux $\Phi^t(V) \in Q= [-1,1]^n $ defines an orbit from the initial $V \in B^+$ that intersects the orbit from the initial state $W \in B^+$. Two different orbits of the flux
do not intersect. Then, the two orbits are the same. If necessary changing the roles of $V$ and $W$, we deduce that  $$\Phi^{t_0}(V) = W \mbox{ for some } t_0 \geq 0$$

 $B^+ \subset B $, so $V$ has at least one component $V_k = 0$ and all of them not negative.
But $\Phi_j^t $ is the strictly increasing in time solution of the differential equation $d\Phi_j^t/dt = \gamma_j (\Phi^t_j)$, with $\gamma_j >0$ for all $j $.

We deduce that if $V \in \widehat B^+$, and if  $\Phi^{t_0}(V) = W$ for $t_0 >0$, then $W_j >0 \; \forall \, j$, and therefore $W \not \in B^+$.

As we know that $W \in B^+$ and $t_0 \geq 0$, we conclude that $t_0 = 0$, and then $W = \Phi^{t_0}(V) = \Phi^0(V) = V. \; \; \; \; \Box$

  \begin{definition}
  \label{definicionSeparationProperty} {\bf The Separation Property.} \em
  We say that $F$ verifies \em the separation property \em if $$f_i(B^+_i) \cap f_j(B^+_j) = \emptyset\; \; \; \forall i \neq j$$ where $\{B^+_i\}, \; i = 1, 2, \ldots, n$, are the continuity pieces of $F$ in the positive Poincar\'{e} section $B^+$, as defined in \ref{definicionModeloMatematico}.7., and $f_i$ is the continuous expression of $F|_{B_i^+}$ according to the formulae (\ref{formulaReturnMapF3}).

  Note that  $B^+_i$ is compact for all $i$, and $F$ is  continuous in each $B^+_i$. Therefore the image $F(B_i^+)$ is a compact set. Then, the separation property implies that \em there exists a minimum positive distance $\alpha >0$ between the images by $F$ of two different continuities pieces. \em
  \end{definition}

\begin{theorem}
\label{teoremaSparationProperty}
The Poincar\'{e} map $F$ verifies the separation property.
\end{theorem}

{\em Proof: } Take $B^+_i$ and $B^+_j$ with $i \neq j$. The formulae (\ref{formulaReturnMapF3}) of the Poincar\'{e} map $F|_{B^+}: B^+ \subset B \mapsto B$ and the Theorem \ref{teoremaBmasinvariante} imply that
$$\forall \, \, V \in B_i^+: \; \; (f_i)_i(V) = 0 , \; \; \; \; (f_i)_j(V) >0 \; \; \forall \, j \neq i $$
$$\forall \, \, W \in B_j^+: \; \; (f_j)_j(V) = 0 , \; \; \; \; (f_j)_i(V) >0 \; \; \forall \, i \neq j $$
Then $f_i(B_i^+) \bigcap f_j(B_j^+) = \emptyset. \; \; \; \Box$
\begin{remark}
 \label{remarkGloballyInjectiveness} {\bf Global injectiveness of the Poincar\'{e} map.}

 \em From Theorems \ref{teoremaFinjective} and \ref{teoremaSparationProperty} it is deduced that the Poincar\'{e} map $F$ is globally injective in $B^+$. In fact, if $V \neq W$ are in the same continuity piece $B_i^+$, then $f_i(V) \neq f_i(W)$ because $F|_{B_i^+} = f_i$ is injective. And if $V \neq W$ respectively belong to two different continuity pieces $B_i^+$ and $B_j^+$ for $i \neq j$, then $f_i(V) \neq f_j(W)$ because $f_i(B_i^+) \cap f_j(B_j^+) = \emptyset$, due to the separation property. We deduce that if $V \neq W$ then $F(V) \bigcap F(W) = \emptyset$, where the image set $F(V)$ of a point is defined in the Remark \ref{definicionModeloMatematico}.4.
 \end{remark}
\begin{theorem}
\label{teoremaContraccion} {\bf Local contractiveness.}

The Poincar\'{e} map $F|_{B^+_i}$ is uniformly contractive, but not infinitely contractive, in each of its continuity pieces $B^+_i$. Precisely, there exist two constant real numbers $0 < \sigma < \lambda < 1$ and a distance $\dist$ in the positive Poincar\'{e} section $B^+ = \bigcup B_i^+$, such that, for all $i = 1, 2, \ldots, n$:
$$\sigma \dist (V,W) \leq \dist (f_i(V), f_i(W)) \leq \lambda \dist (V,W) \; \; \; \forall \, V, W \in B_i^+$$ where $f_i: B_i^+ \mapsto B^+$ is the continuous restriction of $F$ to $B_i^+$, according with formulae (\ref{formulaReturnMapF3}).
\end{theorem}
{\bf Remark: } The distance $\dist$ of Theorem \ref{teoremaContraccion} induces the same topology in $B^+$ as homemorphic to a compact ball of $\mathbb{R}^{n-1} \subset \mathbb{R}^n$. In fact, along the proof of the Theorem \ref{teoremaContraccion} we will construct a linear projection $\pi : \mathbb{R}^n \mapsto \mathbb{R}^{n-1} $ and a diffeomorphism $\xi: \mathbb{R}^n \mapsto \mathbb{R}^n $ of $C^1$ class, such that:
$$\dist(V, V + dV) = \| \pi ( d \xi dV)\|  \mbox{ where }
\| \cdot \| \mbox{ is a norm in }\mathbb{R}^n. $$

\vspace{.3cm}

{\em Proof of the Theorem  } \ref{teoremaContraccion}: The continuity piece $B_i^+$is fixed. For simplicity of the notation, along this proof we will use simply $F$ to denote $f_i$.

The existence of the distance $\dist$ and the contraction rate $\lambda$ is proved in the Theorem 3 of \cite{Yo}. For a seek of completeness we include here some pieces of the proof of \cite{Yo}, adding to them  the existence of the bound contraction rate $0<\sigma <1$, $\sigma < \lambda$.

Due to the Tubular Flux Theorem there exists a ${\cal C}^1$   diffeomorphism which is a spatial change of variables $\xi: V \mapsto  \breve{V} $ from $Q \subset \mathbb{R}^n$ onto $\breve Q \subset \mathbb{R}^n$, such that $\xi|_{B^+} = id$ and the solutions of the differential equation (\ref{ecuacionDiferencial}) in $Q$ verify $$d \breve V/dt = \vec a $$ in $\breve Q$, where $\vec a \in \mathbb{R}^n$ is a constant vector with positive components. It verifies: $$\xi (\phi^t(V)) = \xi (V) + \vec a \cdot t, \; \;
d\xi \cdot \gamma (V) = \vec a  \; \; \forall \, V \in Q$$

   Define in $\mathbb{R}^n$ the ortogonal projection $\pi$
 onto the $(n-1)$-dimensional subspace
  $$a_1 \breve V_1 + a_2 \breve V_2 + \ldots + a_n \breve V_n = 0$$ The flux of the differential equation (\ref{ecuacionDiferencial}), after the change $\xi $ of variables in the space,  is ortogonal to that subspace, and is transversal to $\xi (B^+) = \breve B^+ = B^+ $ (recall that $\xi|_{B^+} $ is the identity map).

  Consider any real function $g: \mathbb{R}^{(n-1)} \mapsto \mathbb{R}$:  \begin{eqnarray} \forall\; \breve{V},  \breve{V} + d\breve{V} \in \mathbb{R}^n:\; \; \; \; \;  \pi(d\breve{V})=  \pi (d\breve{V} +   g (\breve{V}) \cdot \vec a) . && \nonumber
  \\ \forall \, V, \, V + dV, \; \; U  \in \breve B_k, \mbox{ define } \; \; \;  \dist(V, V + dV) = \| \pi ( d \xi \,  dV)\| && \nonumber \\ \dist( V,U ) = \int _{0}^1 \|\pi ( d \xi _{V  + t (U-V)} \cdot (U-V)  \| \, dt && \label{definicionDistanciaContractiva} \end{eqnarray}
   It is left to prove that $f_i: B_i^+ \mapsto B^+$ is contractive with this distance.

   Let us apply $f_i$  to $V$ and $ V + dV $ in $\in B_i^+ $. We use the equalities (\ref{formulaReturnMapF3}).

   We shall use the Liouville derivation formula of the flux of the differential equation  respect to its initial state: \begin{eqnarray} &{d \Phi_j ^t }/{ d V_j} = \exp \left ( {\int _0 ^t } \gamma'_j  \, (\Phi_j^s (V_j)) \, ds) \right ) & \nonumber \end{eqnarray}
Define: $$ -\alpha = \max_{j} \, \max _{V_j \in [-1,1]} \; \; { \gamma'_j (V_j)}\; \; < 0 , \; \; \; \; \; \;  -\alpha ^* = \min_{j} \, \min _{V_j \in [-1,1]} \; \; { \gamma'_j (V_j)}\; \; < 0$$
   Use the Lemma \ref{lemaTiempoMinimo} to bound uniformly above zero the inter-spike intervals $\overline t(V)$: $$0 < T \leq \overline t(V)$$ Recall that $\overline t(V)$ is the solution of the $C^1$ implicit equation $\Phi ^{\overline t (V_i) } (V_i) = 1$. Then $\overline t (V) $ is a continuous real function of $V \in B^+$, and $B^+$ is a compact set. So, $\overline t (V)$ is also upper bounded by a constant: $$ \overline t(V) \leq T^*$$

   Derive the formulae (\ref{formulaReturnMapF3}) to obtain:
  \begin{eqnarray}&F(V) - F(V + dV )  = d F\cdot dV =  \left [ ({\partial F_j }/{\partial V_j}) dV_j + ({\partial F_j }/{\partial V_i}) d V_i   \right ] _{1 \leq j \leq n}& \nonumber \\ && \nonumber \\ &
  {\partial F_j }/{\partial V_j} \,  = \left . ({d \Phi_j^{t}(V_j) }/{d V_j} ) \right |_{t = \overline t (V)}  = \exp \left ( \int _0 ^{\overline t(V)}  \gamma'_j (\Phi_j^s (V_j)) \, ds) \right )  \in [e^{-\alpha ^* T^*},  e^{- \alpha  T}] \nonumber &  \\ &  & \label{formulaNombre}\\ & {\partial F_j }/{\partial V_i} \,  =   \left . ({d \Phi_j^{t}(V_j) }/{dt}) \right |_{t = \overline t (V)}  \cdot ({d t_i (V_i)}/{dV_i})=  g(V) \cdot \gamma_j(\Phi_j^{\overline t (V)}(V_j)) &  \nonumber \end{eqnarray}
  where $g(V) = {d t_i (V_i)}/{dV_i}$ is the real function obtained deriving respect to $V_i$ the implicit equation given in (\ref{ecuacionTiempoSpike}): $1 = \Phi_i^{t_i (V_i)}(V_i) $. Call $\vec e_j$ to the $j-$th. vector of the canonic base in $\mathbb{R}^n$ and join all the results above:

   \begin{eqnarray} & \pi \cdot  d \xi \, \left ( F(V + dV) - F(V  ) \right ) = \pi \cdot  d \xi d F \cdot dV =  & \nonumber \\ & = \pi \cdot d \xi  \left (\sum_{j= 1}^n ({\partial F_j }/{ \partial V_j})  \cdot dV_j \vec e_j \right )   + \pi \cdot d \xi  (g(V) \cdot  \gamma (\Phi^{\overline t(V)}(V) ) = & \nonumber \\ & = \pi \cdot d \xi  \left (\sum_{j= 1}^n({\partial F_j }/{ \partial V_j})  \cdot dV_j \, \vec e_j \right ) +  g(V) \pi \cdot d \xi  \cdot \gamma(\Phi^{\overline t(V)}(V)) )  = & \nonumber \\ & = \pi \cdot d \xi \left (\sum_{j= 1}^n({\partial F_j }/{ \partial V_j})  \cdot dV_j \, \vec e_j \right ) +  \pi ( g(V) \cdot \vec a )= \pi \cdot d \xi  \left (\sum_{j= 1}^n({\partial F_j }/{ \partial V_j} ) \cdot dV_j \, \vec e_j \right ) & \label{ecuacionDiferencialDeRho}\end{eqnarray}  We define the numbers $ \sigma$ and $ \lambda $:  $0 <\sigma = e^{-\alpha^* T^*} <  e^{-\alpha T} = \lambda <1$ and observe from the computations in (\ref{formulaNombre}) that:
  $$0 < \sigma = e^{-\alpha^* T^*} \leq \partial F_j / \partial V_j \leq   e^{-\alpha T} = \lambda < 1      $$

\noindent Applying the definition of the differential distance $\dist $ in (\ref{definicionDistanciaContractiva}),  and the equality (\ref{ecuacionDiferencialDeRho}), we obtain:
\begin{eqnarray}
&\dist(F(V), F(V + dV) ) = \| \pi (d \xi \cdot  d F \cdot dV) \| \leq   \lambda  \, \|\pi d \xi \cdot dV  \| = \lambda \, \dist(V, V + dV) = \lambda \, \| \pi ( d \xi \, dV)\| & \nonumber \\
&\dist(F(V), F(V + dV) ) = \| \pi (d \xi \cdot  d F \cdot dV) \| \geq   \sigma \, \|\pi d \xi \cdot dV  \| = \sigma \, \dist(V, V + dV) = \sigma \, \| \pi ( d \xi \,  dV)\| & \nonumber \end{eqnarray}
Integrating by formula (\ref{definicionDistanciaContractiva}) we conclude:
$$ \sigma \, \dist (V, U) \leq \dist ( F(V) , F (U) ) \leq \lambda \, \dist (V, U) \; \; \Box$$

\begin{remark} {\bf Local homeomorphic property of the Poincar\'{e} map. }

{ Each continuity piece $f_i$ of the Poincar\'{e} map in $B_i^+$ is an homeomorphism onto its image.}
\label{remarkhomeomorfismo}

\em
It is  an immediate consequence of Theorem \ref{teoremaContraccion} and the global injectiveness of $F$. Even more, the continuous restriction $f_i = F|_{B_i^+}$,  is Lipschitz with constant $\lambda < 1$ and its inverse (defined from $f_i(B_i) \mapsto B_i$) is also Lipschitz with constant $1/\sigma >1$. Then $f_i$ is an \em homeomorphism \em onto its image. $\Box$

We note that the same result can be obtained without using the Theorem \ref{teoremaContraccion}. Due to Theorem \ref{teoremaFinjective}, $f_i$ is injective, and due to the formulae (\ref{formulaReturnMapF3}), $f_i$ is continuous. Due to the Theorem of the invariance of the Domain, any continuous and injective function from a ball in $\mathbb{R}^{n-1}$ to $\mathbb{R}^{n-1}$ is an homeomorphism onto its image.
\end{remark}
\vspace{.5cm}

In the following corollary we resume all the conclusions of this section:

\begin{corollary}
\label{corolario} If the network  of $n$ inhibitory neurons verifies the assumptions of the physical model, \em evolving with real time $t$ in the phase space $Q \subset \mathbb{R}^n$ as stated in (\ref{ecuacionDiferencial}), (\ref{ecuacionSinapsis1}), (\ref{ecuacionSinapsis2}), (\ref{formulaEpsilo0GrafoCompleto}), (\ref{ecuacionesLargeDissipativeness1}), (\ref{ecuacionesLargeDissipativeness2}) and (\ref{ecuacionesLargeDissipativeness3}), \em then there exists a Poincar\'{e} section $B^+$ \em   that is homeomorphic to a $(n-1)$-dimensional compact ball, \em and a return map $F: B^+ \mapsto B^+$, \em which has the following properties:

  {\bf a)} \em $F$ is piecewise continuous. \em Precisely: there exist a finite partition $   \{B_i^+\}_{1 \leq i \leq n}$ of the Poincar\'{e} section $B^+$, formed by compact sets $B_i^+$ homeormorphic to compact balls of ${\mathbb{R}^{n-1}}$, with pairwise disjoint interiors, and there exist $n$ continuous maps $f_i : B_i^+ \mapsto B^+$ being $F(V) = \{ f_i(V) : \; i \mbox{ such that } V \in B_i^+\} $ for all $V \in B^+$.

   As a consequence $F$  is univoquely defined as $f_i$ in the interior  of its continuity piece $B_i$, and  multidefined as $f_{i}, f_{j} $ in $B_i \cap B_j$, if $i \neq j$, $f_i \neq f_j$.

  {\bf b)} \em $F$ is locally uniformly contractive and not infinitely contractive, \em i.e. for some metric $\dist$ in $B^+$ the exist  constants $0 < \sigma < \lambda < 1$ such that for all $1 \leq i \leq n$: $\sigma \dist (V,W) \leq \dist (f_i(V), f_i(W)) \leq \lambda \dist (V, W) \; \; \forall \, V, W \in B_i^+$.

  {\bf c)} \em $F$ has the separation property, \em i.e. $f_i(B_i^+) \cap f_j(B_j^+) = \emptyset$ if $i \neq j$.  Therefore, there exists $0 < \alpha = \min _{i \neq j} \dist (f_i(B_i^+), f_j(B_j^+))$.

  Note that from b) and c), it is deduced that $F$ is globally injective in $B^+$, as proved in Remark \ref{remarkGloballyInjectiveness}. Also from b) it is deduced that   $f_i : B_i^+ \mapsto f_i (B_i^+) \subset B^+$ is an homeomorphism onto its image.

\end{corollary}

\begin{Paragraph}
{\bf Comments about the mathematical model.}
\end{Paragraph}
Due to Corollary \ref{corolario}, all the general results that we will prove for abstract piecewise continuous maps $F$ verifying a), b), c), are applicable to the networks of inhibitory neurons in the assumptions of the physical model stated in \ref{definicionModeloInicial}. Nevertheless the reciprocal of the Corollary \ref{corolario} does not hold. Given a map $F$ verifying a), b), c) there does not necessarily exist a  network of inhibitory neurons in the hypothesis of the physical model stated in \ref{definicionModeloInicial} for which $F$ is its first return Poincar\'{e} map.

Nevertheless we can wide our scenario of possible  inhibitory neuronal networks models. In fact, the properties a), b) c) are open (in the uniform ${\cal C}^0 + {Lipschitz}$ topology of the finite family of maps $f_i$). Thus they are not only verified by systems for which the differential equations (\ref{ecuacionDiferencial}) are independent in the $n$ variables $V_i$, but also if the system is of the form $dV/dt = \gamma^* (V)$, where $\gamma^* : \mathbb{R}^n \mapsto \mathbb{R}^n$ is a ${\cal C}^1$ vector field, near enough  the given  $\gamma = ( \gamma= \gamma_1, \gamma_2, \ldots, \gamma _n)$, even if $\gamma^*$ does not verify all the hypothesis stated in  (\ref{ecuacionDiferencial}).

Also the matrix $(H_{i,j})_{i,j}$ of synaptic interactions in the network can be substituted for any matrix $(H^*_{i,j})_{i,j} (V)$, not necessarily constant, but functions near the constant matrix $(H_{i,j})_{i,j}$ and so, still verifying the assumptions (\ref{formulaEpsilo0GrafoCompleto}), (\ref{ecuacionesLargeDissipativeness1}), (\ref{ecuacionesLargeDissipativeness2}), (\ref{ecuacionesLargeDissipativeness3}). Therefore, without changing the synaptical rules in equations
(\ref{ecuacionSinapsis1}) and (\ref{ecuacionSinapsis2}), but allowing the synaptic interactions slightly depend of the postsynaptic potentials, we will obtain a Poincar\'{e} map $F$ still verifying the thesis a), b), c)  of the Corollary \ref{corolario}.

Besides, as observed in the subsection \ref{comments}, the physical model  includes  looser hypothesis than those specified in \ref{definicionModeloInicial}, modulus any differentiable change of the variables of the system. So, also in those  models the properties a), b), c)  are verified by an open family of systems.

Finally, the properties a) b) c) of the Corollary \ref{corolario} are verified by many other models, in which the interspike regime is stated as a dynamical system depending continuously on time $t$ and on the initial state $V$, but not necessarily as regular as to verify a differential equation.  The dynamics of the potential $V_i$ in the inter-spike interval may be  given by a flux $\Phi_i ^t (V_i )$  defined continuously in time $t$,  strictly increasing on $t$, continuous  but not necessarily differentiable respect to $t$ nor to the initial state. But not all such general models are in the aim of this work. They must be posed some hypothesis, to get the properties a) b) and c) of the  Poincar\'{e} section and its return map $F$.

   The arguments above justify to wide  the abstract mathematical model of a network of $n$ inhibitory neurons, according to the following definition:

\begin{definition}
\label{definicionModeloAbstracto} {\bf The  Abstract Mathematical Model.} \em
We say that a map $F: B^+ \mapsto B^+$, in a set $B^+$ homeomorphic to a compact ball of $\mathbb{R}^{n-1}$, \em models a generalized network of $n $ inhibitory neurons \em if it verifies the statements a), b), c)  of the Corollary \ref{corolario}.
\end{definition}

\section{The abstract dynamical system.} \label{definicionesabstractas}

 Let $B \subset {\mathbb{R}^{n}}$  be a compact set,  homeomorphic  to a compact ball of $\mathbb{R}^{n-1}$. In particular $B_i$ is connected.

\begin{definition}
\label{particion}

\em
A \em finite partition \em of $B$ is a finite collection
$\{B_i\}_{1\leq i \leq m}$  of
compact non empty sets $B_i$ of $B$, homeomorphic to compact balls of ${\mathbb{R}^{n-1}}$,
 such
that
 $ \bigcup _{1\leq i \leq m}B_i = B$  and  $\mbox{ int }B_i \, \cap
 \mbox{ int }B_j \, = \, \emptyset $, for $ i \neq j $.

 Denote $S =  \bigcup _{i \neq j } B_i \cap B_j$, and call $S$
the \em separation line, \em  or \em line of discontinuities \em, although it is not a line in the usual sense, but the union
of the  topological frontiers of $B_i$, each one homeomorphic to some $(n-2)$-dimensional manifold.
\end{definition}

\begin{definition} \em \label{piecewisecontinuous}

Given a finite partition $\{B_i\}_{1\leq i \leq m} $ of $B$,
we call $F$ a   \em  piecewise
 continuous map on $(B, {\cal P})$ with the separation property \em
  if $F$ is a finite family  $
 F = \{f_i\}_{1\leq i \leq m}$
 of  \em homeomorphisms \em $f_i: B_i \mapsto f_i(B_i) \subset B $, such that $f_i (B_i) \cap
 f_j (B_j) = \emptyset$ if $i \neq j$. We note that $F$ is multi-defined in the separation line $S$.

  Each $B_i$ shall be called a \em continuity piece \em of $F$.

  \end{definition}

   \begin{remark}
   \label{remarkFinversa}
    \em  A piecewise continuous map $F$ with the separation property is globally injective because it is an homeomorphism in each continuity piece and two different continuities pieces have disjoint images. Therefore $F^{-1}$ exists, uniquely defined in each point of  $F(B) = \bigcup_i f_i(B_i)$. In fact:

   For any point $x \in \bigcup_i f_i(B_i)$, its backward first iterate is uniquely defined as  $F^{-1}(x) = f_i^{-1}(x)$, where $i$ is the unique index value such that $x \in f_i(B_i)$.

   Nevertheless $F^{-1}$ is  not necessarily injective because  $F$ is multidefined in $S = \bigcup_{i \neq j} (B_i \cap B_j)$.

    $F^{-1}$ is continuous in $F(B)$, because $F^{-1}|_{f_i(B_i)} = f_i^{-1}$ and $f_i$ is an homeomorphism due to the Definition \ref{piecewisecontinuous}.
   \end{remark}

\begin{definition} \em  \label{definicioncontraccionlocal}

We say that $F$ is  uniformly \em locally contractive  \em if there exists a constant $0 < \lambda < 1$, called an \em uniform contraction rate for $F$, \em  such that
$\dist (f_i(x),f_i(y)) \leq \lambda \dist (x,y)$, for all $x$ and $y$ in the same $B_i$,
for all $1 \leq i \leq m $ .

 \end{definition}

Given a point $x \in B$,     its image set  is $F(x) = \{f_i(x):   x \in B_i \}$.
If $H \subset B$, its image set is $F(H) = \bigcup _{x \in H} F(x)$.
We have that $B \supset F(B) \supset \ldots F^k(B) \supset \ldots $.

The second iterate of the point $x \in B$ is the set $F^2(x) = F (F(x))$. Analogously is defined the $j-$th. iterate as the set $F^{j}(x)$ for any $j \geq 1$. We convene to define $F^0(x) = x $ and $F^0(H) = H$.

\begin{definition} \em \label{definicionAtomo}
For any natural number $k \geq 1$, we call \em atom of generation $k$ \em
to
 $$f_{i_k}\circ \ldots \circ f_{i_2} \circ
 f_{i_1} (B_{\mathbb{I}}) $$ where $\mathbb{I} = (i_1, i_2, \ldots, i_k) \in \{1,2,\ldots, m\}^k$
 and $B_{\mathbb{I}}$ is  the subset of $B_{i_1}$ where the  composed function above is defined. (If $B_{\mathbb{I}}$ were an empty set, then the atom is empty.)  Abusing of the notation we write the atom as:
 $$f_{i_k}\circ \ldots \circ f_{i_2} \circ
 f_{i_1} (B_{i_1})$$
 \end{definition}

 We note that each atom of generation $k$ is a compact,
  not necessarily connected set, whose
 diameter is smaller than $\lambda ^k \mbox{diam} B$.

 The set $F^k(B) $ is a compact set, formed by the union of all the
 not empty atoms
of generation $k$.

There are at most $m^k$ and at least $m$ not empty atoms of generation $k$, where $m$ is the number of continuity pieces of $F$.

\begin{definition} \label{limitset}
\em

Given $x_0 \in B$, a future orbit $o^+(x_0)$ is a sequence of points $\{x_i\}_{i \geq 0}$, starting in $x_0$, such that $x_{i+1} \in F(x_i) \; \; \forall \, i \geq 0$. Due to the multi-definition of $F$ in the separation line $S$, the points of $S$ and those that eventually fall in $S$ may have more than one future orbit.

A point $y$ is in the limit set $L^+(o^+(x_0))$ of a future orbit of $x_0$ if there exists $k_j \rightarrow +\infty$ such that $x_{k_j} \rightarrow y$.

The limit set $L^+(x_0)$ is the union of the limit sets of all its future orbits.

  The \em limit set \em $L^+(B)$ of the map $F$, also denoted as $L^+(F)$, is
 the union of the limit sets of all the points $x \in B$. \em
 \end{definition}

 \begin{remark} \em
 \label{remarkConjLimiteinvariante} Due to the compactness of the space $B$ the limit set $L^+(o^+(x_0))$ of any future orbit, \em is not empty. \em

 Also, it is standard to prove that \em $L^+(o^+(x_0))$ is compact \em (because it is closed in the compact space $B$). Nevertheless $L^+(x_0)$ may be not compact, if the point $x_0$ has infinitely many different future orbits.

 Finally, we assert that \em $L^+(o^+(x_0))$ is invariant: \em $F^{-1} (\; L^+(o^+(x_0)) \;  ) = L^+(o^+(x_0))$.

 \end{remark}

 {\em Proof:}  Consider $y \in L^+(o^+(x_0))$. We have $y = \lim _{j \rightarrow + \infty} x_{k_j} \in F(B)$ if $k_j \geq 1$.

  $F^{-1}: F(B) \rightarrow B$ is a continuous uniquely defined function in the compact set $F(B)$ (see Remark \ref{remarkFinversa}).  Then $x_{k_j -1} = F^{-1} (x_{k_j}) \rightarrow F^{-1}(y) $, so $F^{-1} (y) \in L^+(o^+(x_0))$ proving that $$F^{-1} (\; L^+(o^+(x_0)) \;  ) \subset L^+(o^+(x_0)) $$

 Let us prove the converse inequality: $F^{-1} (\; L^+(o^+(x_0)) \;  ) \supset L^+(o^+(x_0)) $.

 $F = \{f_i: B_i \mapsto B\}$ is defined and continuous in each of its finite number of pieces $B_i$ that are compact and cover $B$. Then there exists some $i \in\{1, 2, \ldots, n\}$ and a subsequence (that we still call $k_j$), such that $$y = \lim_{j \rightarrow + \infty} x_{k_j} \in B_i,  \; \; \; \forall j \geq 0: \; \; x_{k_j} \in B_i, \; \; x_{k_j + 1} = f_i (x_{k_j}), \; \; f_i(y) = \lim f_i(x_{k_j}) = \lim x_{k_j + 1} $$
 There exists $y_1 = f_i(y) \in F(y)$ such that $y_1 \in L^+(o^+(x_0))$. In other words, $y \in F^{-1}(L^+(o^+(x_0))$. This last assertion was proved for any $y \in L^+(o^+(x_0))$. Therefore $ L^+(o^+(x_0)) \subset F^{-1} (\; L^+(o^+(x_0)) \;  ) $ as wanted. $ \; \Box$
 \begin{definition}

\em
 We say that a point $x$ is \em  periodic of period $p$  \em
 if there exists a first natural number $p\geq 1$   such that $x \in F^{p}(x) $. This is
 equivalent to $x$ be a periodic point in the usual sense, for
  the uniquely defined map $F^{-1}$, i.e.
  $F^{-p}(x)= x$ for some first natural number $p \geq 1$.

  We call the backward orbit
 of $x$ (i.e. $\{ F^{-j}(x), j = 1,\ldots, p\}$), a periodic orbit
  with period $p$.

\end{definition}

   We will prove in Lemma \ref{lemaprevio} that the limit set $L^+(B)$ is contained in the compact, totally
   disconnected set $K_0 = \bigcap _{k\geq1} F^k(B)$. It could be
   a Cantor set. But generically
    $K_0$  shall be the union of a finite number of periodic orbits,
    as  we shall prove in Theorem \ref{Teorema1}.

\begin{definition} \em  \label{finallyperiodic}

  We say that $F$ is \em finally periodic \em with period $p$ if the
  limit set $L^+ (F)$ is the union of \em only \em a finite number of periodic orbits
  with minimum common multiple of their periods equal to $p$.
 In this case we call \em limit cycles \em to the periodic orbits of  $F$.

 We call basin of attraction of each limit
 cycle  $L$ to the set of points $x \in B$ whose  limit set $L^+(x)$ is $L$.

\end{definition}

\vspace{.3cm}

{\bf Topology in the space of piecewise continuous locally contractive maps in $B$.}

\vspace{.3cm}

  Let ${\cal P} = \{B_i\}_{1\leq i \leq m}$ and  ${\cal Q} =
  \{A_i\}_{1\leq i \leq m}$  be   finite partitions (see Definition \ref{particion})
  of the compact region $B$ with the same number $m$ of pieces.

We define the distance between ${\cal P}$ and ${\cal Q}$ as
\begin{equation}
\label{definicionHausdorfDistance}
d({\cal P}, {\cal Q}) = \max_{1 \leq i \leq m}\; \; \;  \mbox{Hdist} (A_i, B_i) \end{equation} where $\mbox{Hdist} (A, B)$ denotes the Hausdorff distance between the two sets $A$ and $B$. i.e.
$$\mbox{Hdist} (A, B) = \max \{\dist (x, B), \dist (y,A), x \in A, y \in B  \}$$

\begin{definition} \em  \label{topologia}
  Let
   $
  F  = \{f_i:B_i \mapsto B\}_{1\leq i \leq m}$
  and $G   =\{g_i:A_i \mapsto B\}_{1\leq i \leq m}$  be    locally contractive
   piecewise continuous maps on $(B, {\cal P})$ and  $(B, {\cal Q})$ respectively.
    Given
   $\epsilon >0 $
  we say that   $G $ is a
  \em $\epsilon$-perturbation
 of $F$  \em if
$$ \max_{1 \leq i \leq m} \left \| \left. (g_i - f_i) \right | _ {\displaystyle {B_i \cap A_i}} \right \|_{\displaystyle {{\cal C\;} ^0}} < \epsilon, \; \; |{ \lambda_F - \lambda_G} | < \epsilon \; \; \; \mbox { and
} \; \; d({\cal P}, {\cal Q} ) < \epsilon $$
where $\lambda _F$ denotes \em the uniform contraction rate \em of $F$ in  its continuity pieces, defined in \ref{definicioncontraccionlocal}, and $\|\cdot \|_{{\cal C}^0}$ denotes the  ${\cal C}^0$ distance in the functional space of continuous functions defined in a \em compact \em set $K$: $$\| (g - f)|_K\|_{{\cal C}^0} = \max_{ x \in K} \dist (g(x), f(x))$$
\end{definition}

\begin{definition} \label{definicionPersistencia} \em
 We say that the limit cycles  of a finally periodic map $F$  (see Definition \ref{finallyperiodic})
  are \em persistent  \em if:

  For all $\epsilon^*>0$ there exists $\epsilon >0$
   such that all $\epsilon$-perturbations $G$
of $F$  are finally periodic  with the same
finite number of limit cycles (periodic orbits) than $F$, and such that each limit  cycle $L_{G}$ of $G$ has the same period and is $\epsilon^*$-near of some limit cycle $L_{F}$ of $F$ (i.e. the Hausdorff distance between $L_G$ and  $L_F$ verifies
$\mbox{Hdist}(L_G, L_F) < \epsilon^*$).

\end{definition}

\begin{definition} \em  \label{generico}

Denote ${\cal S}$ to the \em space of all the systems that are  piecewise continuous with the separation property and locally contractive, \em according with the Definitions \ref{piecewisecontinuous} and \ref{definicioncontraccionlocal}.

We say that a property $\mathbb{P}$  of the systems in ${\cal S}$ (for instance being finally periodic as will be shown in Theorem \ref{Teorema1}) is \em (topologically) generic \em
if it is verified, at least, by an \em open and dense \em subfamily of systems in the functional space ${\cal S}$, with the topology  (in ${\cal S}$) defined in \ref{topologia}.

Precisely, being \em generic \em means:

1) The \em openness \em  condition: For each piecewise continuous map $F$ that verifies the property $\mathbb{P}$ there exist $\epsilon >0$ such that
all $\epsilon$-perturbation of $F$ also verifies $\mathbb{P}$.

2) The \em denseness \em condition:  For each piecewise continuous map $F$ that does not verify the property $\mathbb{P}$, given $\epsilon >0$, arbitrarily small, there exist some $\epsilon$-perturbation $G$ of $F$ such that $G$ verifies the property $\mathbb{P}$.

\vspace{.2cm}

The openness  condition implies that the property $\mathbb{P}$ shall be robust under small perturbations of the system. It is robust under small changes, not only of a finite number of real parameters, but also of the functional parameter that defines the model itself. So   the system should be structurally stable. When this robustness holds, the property $\mathbb{P}$ is still observed  when the   system,  the model itself, does not stay  exactly fixed, but is changed, even  in some unknown fashion, remaining near the original one.

The density condition combined with the openness condition,  means that the only behavior that have chance to be observed under not exact experiments are those that verify the property  $\mathbb{P}$. In fact, if the system did not exhibit the property $\mathbb{P}$, then some arbitrarily small change of it, would lead it to exhibit $\mathbb{P}$ robustly. \em

\end{definition}

The denseness condition implies that if the property $\mathbb{P}$ were generic, then the opposite property   (Non-$\mathbb{P}$) has null interior in the space of ${\cal S}$ of systems, i.e. Non-$\mathbb{P}$ is not robust: some arbitrarily small change in the  system will lead it to exhibit $\mathbb{P}$. That is why we define the following:

\begin{definition} \em  \label{bifurcating}
If the property $\mathbb{P}$ is generic, we say that any system that does not exhibit $\mathbb{P}$ is \em bifurcating, \em  and  Non-$\mathbb{P}$ is a \em not persistent \em property.
\end{definition}

\section{The generic persistent periodic behavior.} \label{seccionperiodico}

\begin{theorem} \label{Teorema1}
Let $F$ be a locally contractive  piecewise continuous
map with the separation property.  Then generically $F$ is finally periodic with persistent limit cycles. \em

\end{theorem}

To prove Theorem \ref{Teorema1} we shall use the following lemma:

 \begin{lemma} \label{lema}
If there exists an integer $k \geq 1$ such that the compact set
 $K =  F ^k (B)$ does not intersect the separation line $S$ of the partition into the continuity pieces of $F$, then
 $F$ is finally periodic and its limit cycles are persistent.

\end{lemma}

{\em Proof: }
 By hypothesis, $\dist (K, S) = d >0$, because $K$ and $S$ are disjoint compact sets. On the other hand $K= F^k(B) = \bigcup _{A \in {\cal A}_k} A$, where ${\cal A}_k$ denotes the family of all the atoms of generation $k$.

 As the  diameter of
 each of the finite number of atoms of generation $k$  is smaller than $diam (B) \, \lambda^k $, it converges to zero when $k \rightarrow +\infty$. Thus, for all  $k$ large enough, it is smaller than $d/2$.

 We assert that each atom $A$ of such generation $k$, is contained in the interior of some continuity piece $B_i$. In fact, fix a point $x \in A$. As the continuities pieces cover the space $B$, there exists some (a priori not necessarily unique) index $i$ such that $x \in B_i$. It is enough to prove that $y \in \    int (B_i)$ for all $y \in A$ (including $x$ itself).

 We argue  in the compact and connected metric space $B$, using known properties of any general compact and connected metric space, for instance  the triangular property, and also the property asserting that the distance of a point $y$ to a set, is the same that the distance of $y$ to the frontier of that set.

  We denote $B_i^c$ to the complement of $B_i$ in $B$, and in the topology relative to $B$  we denote: $\overline {(B_i^c)}$ to the closure of $B_i^c$, i.e the complement of $int(B_i)$, and $\partial B_i$ to the frontier of $B_i$ in $B$, $\partial B_i \subset S$:
 $$\dist (x,y) \leq diam (A) < d/2, \; \; \dist (x, \overline {(B_i^c)}) = \dist (x, \partial B_i) \geq \dist (x, S) \geq d$$ $$ \dist (y, \overline {(B_i^c)}) \geq \dist (x, \overline{(B_i^c)}) - \dist (x,y)  \geq d - d/2 = d/2 >0 $$
 Therefore $y \not \in \overline{(B_i^c)}$ proving the assertion.

  We deduce that given an atom $A \in {\cal A}_k$, there exists and is unique a  natural number $i_0$ such that $A \in int (B_{i_0})$. Therefore $F(A)$ is a single atom of generation $k+1$.

  From the definition of atom in \ref{definicionAtomo}, we obtain  that any atom of generation larger than $k$ is contained in an atom of generation $k$. But each atom of generation $k$ is in the interior of a piece of continuity of the partition $\{B_i\}$.

  We deduce that there exists a sequence of natural numbers $\{i_h\} _{h \geq 0}$, such that \begin{equation}
  \label{equationImagenesDeAtomo}
  A \in int B_{i_0 }, \; \; F(A) = f_{i_0}(A) \subset int B_{i_1}, \; \; \; F^2(A) =  f_{i_1} \circ f_{i_0} (A) \subset int B_{i_3}, \ldots\end{equation} and the successive images of the atom $A$ of generation $k$, are  single
   atoms of generation $k+1, k+2, \ldots, k+h, \ldots $.  Therefore, the successive images of the atom $A$, in the sequence (\ref{equationImagenesDeAtomo}), are contained in a sequence of atoms: $A= A_0, A_1, A_2, \ldots, A_h, \ldots, $ all of generation $k$.

   The same property holds for any of these atoms of generation $k$, and each of them is contained in the interior of a continuity piece of $F$, so $F$ is uniquely defined there and we have:
   \begin{equation}
   \label{ecuacionChain2}
   A = A_0 \subset int B_{i_0}, \; \; F(A_0) \subset A_1 \subset int B_{i_1} , \; \;F^2(A_0) \subset F(A_1) \subset A_2  \subset int B_{i_2} , \ldots, \end{equation}

   The family of atoms of generation $k$ is finite, so we conclude that there exists two first natural numbers $0 \leq h < h+p$ such that $F^p (A_h) \subset A_h$.

   Note that,  $F^p (A_h) $ is uniquely defined as $ f_{i_{h+p}} \circ f_{i_{h+p -1}} \circ \ldots \circ f_{i_{h}}$, because we are considering sets contained in the interior of the continuity pieces of $F$.

    Due to the uniform contractiveness of $f_i$ in each of its continuities piece, $F^p : A_h \mapsto A_h$, is uniformly contractive. The Brower Theorem of the Fixed Point states that in a complete metric space, any uniformly contractive map from a compact set to itself, has an unique fixed point, and  all the orbits in the set converge to this fixed point in the future. Therefore, there exists in $A_h$ a periodic point by $F$ of period $p \geq 1$, and all the orbits with initial states in $A_h$ have the periodic orbit $L$ of $p$, as their  limit set.

    By construction $A_h$ was the image of $A$ by an iterate $F^h$ uniquely defined. So we conclude that the limit set of all the points in the atom $A$ is $L$.

    The construction above can be done starting with any initial atom $A$ of generation $k$. And they are a finite family. We conclude that there exists one and at most a finite number of periodic limit cycles, attracting all the orbits of $\bigcup_{A \in {\cal A}_k} A = F^k(B)$.

    The last assertion implies that the limit set of $B$ is formed by that finite family of periodic limit cycles.

    Finally it is left to prove that the limit cycles are persistent according to the definition \ref{definicionPersistencia}.

    The condition of the hypothesis of this lemma is open in the topology defined in \ref{topologia}, because $K$ and $S$ are compact and at positive distance.

    We assert that the itinerary of each of the atoms $A$ of generation $k$, for $k$ fixed and large enough,   remains unchanged when substituting $F$ by $G$, being $G$ a $\epsilon$-perturbation of $F$ for $\epsilon >0$ small enough.

    In fact $A$  and $\widehat A$ (and also all the other atoms of generation $k$, with $k$ fixed) are contained in the images by $F^k$  or by $G^k$ respectively, of some of their one-to-one corresponding continuity pieces. With $k$ fixed, if $\epsilon>0$ is sufficiently small, they remain at  distance larger than the number  $2d/3>0$ from the  separation lines of $F$ and of $G$ respectively, being $d = \dist (K, S)$.

    Thus, if the generation $k$ is chosen so the atoms  have diameter smaller than $d/3$, repeating the argument at the beginning of this proof we show that $A$ is  in the interior of some continuity piece of $F$, and $\widehat A$ is in  the interior  of the respective correspondent continuity piece of $G$.

    On the other hand, the future iterates of any atom of generation $k$ by $F$, and also by $G$, are contained in the atoms of generation $k$. Therefore the  images  of an atom $A$ or $\widehat A$ of generation $k$, by  all the future iterates of $F$ or of $G$ respectively, are in the interior of their respective one-to-one correspondent continuity pieces. Then the itineraries are the same as we asserted.

    As a consequence,  the indexes $i_0, i_1, i_2, \ldots$ in the finite chain of atoms denoted in (\ref{equationImagenesDeAtomo}) and (\ref{ecuacionChain2}),  remain unchanged, and therefore we deduce the following statement:

     \textbf{A}: \em  The number of periodic orbits in the atoms of generation $k$,  and their periods,  remain unchanged, when substituting $F$ by any $\epsilon$-perturbation $G$, if $\epsilon >0$ is sufficiently small. \em

    It is standard to prove by induction on $k \geq 1$ that for any $\epsilon$-perturbation $G$ of $F$, such that $\lambda + \epsilon = \widehat \lambda  < 1$, each atom $\widehat A$ of generation $k$ for $G$, is at distance smaller than $\sum_{j= 0}^k 2 \epsilon \, \widehat \lambda  ^j < 2 \epsilon/ (1 - \widehat \lambda)  = \epsilon ^* >0$ of the respective atom $A$ for $F$ with the same itinerary.

    Therefore we deduce the following statement:

     \textbf{B}: \em Any periodic point found in an atom $\widehat A  $ of generation $k$ for $G$, is at distance smaller than $\epsilon ^*$ than the respective periodic point found in the correspondent atom $A$ for $F$ with the same itinerary.  \em

    The statements A and B imply that the limit cycles are
    persistent according to Definition \ref{definicionPersistencia}.$\;\; \Box$

 \begin{remark} \em
In the proof of Lemma \ref{lema} we did not use the separation property $f_i(B_i) \cap f_j(B_j) = \emptyset\; \; \forall i \neq j $. At the end of the proof of Lemma \ref{lema} we obtained that  the piecewise continuous and locally contractive systems verifying the thesis of the Lemma \ref{lema}, \ even if they do not have the separation property, contain an open family of systems  in the topology defined in \ref{topologia}. Then:

\em In the space of all the piecewise continuous and locally contractive systems (even if they do not have the separation property), those whose limit set is formed by a finite number of persistent limit cycles form  an open family. \em

 Nevertheless, to prove the genericity of the periodic persistent behavior, we  need to prove that the family of periodic maps is dense in the space of systems. In the following proof, to obtain the density we shall restrict to the space of systems ${\cal S}$   that  verify the separation property.

 \end{remark}

 \begin{remark}
\label{remarkperiodosgrandes} \em
From the proof of Lemma \ref{lema},  the first integer $k \geq 1$ such that $F^{k} (B) \bigcap S = \emptyset$ may be very large, and so the period $p$ may be very large.

In fact, if the system has $n \approx 10^{12}$ neurons, and if no  neuron becomes dead, i.e. it does not eventually remain forever under the threshold level without giving spikes, then the periodic sequences $ i_1, \ldots, i_p$, defined as the itinerary of the periodic limit cycles, have inside the period $p$, at least once each of all the indexes $i \in \{1, 2, \ldots, n\}$. Then $p \geq n \approx 10^{12}$.

As we have shown in the proof of the Lemma \ref{lemaTiempoMinimo}, there exists a minimum time $T >0$ between two consequent spikes. Suppose for instance that $T \approx 10 \; [ms]$ and $n \approx 10^{12}$. The lasting time of the periodic sequence  could be approximately $ 10^{-3} \times 10^{12} [s] = 10^9 [s] \geq 31$ years. So, if most of the neurons did not become dead, the observation of the theoretical periodic behavior of the inhibitory system in the future, could not be practical during a reasonable time of experimentation, and only the irregularities inside the period could be registered, showing the system as virtually chaotic.
\end{remark}

\vspace{.5cm}

{\bf \em Proof of Theorem \ref{Teorema1}. } Due to Lemma \ref{lema} the existence of a finite number of limit cycles attracting all the orbits of the space is verified at least for those systems in the hypothesis of \ref{lema}. At the end of the proof of Lemma \ref{lema} we showed that its hypothesis is an open condition. To prove its genericity it is enough to prove now that the hypothesis of Lemma \ref{lema} is also a dense condition in the space of piecewise continuous contractive maps with the separation property.

Take $F$ being  not
  finally periodic. We shall prove that, for all $\epsilon >0$ there exists a $\epsilon-$ perturbation $G$ of $F$ that verifies the hypothesis of Lemma \ref{lema}, and thus $G$ is finally periodic with persistent limit cycles.

Let be given an arbitrarily  small $\epsilon >0$.

The contractive homeomorphisms $f_i$ of the finite family
  $F = \{f_i: B_i \mapsto B\}_i $, with contraction rate $0 < \lambda <1$,
  can be $C^0$  extended to $F_{\epsilon}= \{f_{i, \epsilon}: U_i \mapsto B\}_i$, where $f_{i, \epsilon}$ is an homeomorphism onto its image defined in compact neighborhoods $U_i \supset B_i$ in $B$ (i.e. $B_i = \overline B_i \subset int (U_i) \subset U_i = \overline U_i \subset B$ where the closures and interiors of the sets are taken in the relative topology of $B \subset \mathbb{R}^{n}$), such that $f_{i, \epsilon}|_{B_i} = f_i$, and such that $f_{i, \epsilon}$ is still contractive,  with a  contraction rate $0 < \lambda ' < 1 $ such that $|\lambda -\lambda'| < \epsilon$.

 The extended
 map $F_{\epsilon } = \{f_{i, \epsilon}: U_i \mapsto B\}_i $, is now multidefined on $\bigcup _{i \neq j } U_i \cap U_j \supset S$.
  The separation property is an open condition, thus the extension $F_{\epsilon}$ still verifies $f_{i, \epsilon} (U_i ) \cap f_{j, \epsilon} (U_j) = \emptyset $ for all $i \neq j $,  if the neighborhoods $U_i$ and $U_j$ are chosen at a sufficiently small Hausdorff distance from their respective pieces $B_i$ and $B_j$, and  $\epsilon>0$ is small enough.

 Call $\epsilon_1>0$ to a positive real number  smaller or equal than $\epsilon$, and also smaller or equal than
  the distance from
 $B_i$ to the complement of $U_i$, for all
 $i = 1,2, \ldots m$. Precisely $$0 < \epsilon_1   = \min \{\epsilon, \; \; \min_{1 \leq i \leq m} \dist (B_i, U_i^c)  \}$$

 Consider the compact sets:
$$K^+ = \bigcap_{k \geq 1}\; \; \;  \bigcup _{(i_1, \ldots, i_k) \in \{1,2 \ldots m\}^k}
f_{i_k, \epsilon} \circ \ldots \circ f_{i_1, \epsilon}(U_{i_1}) \;\;\; \supset \;\;
K$$
$$K = \bigcap_{k \geq 1}\; \; \;  \bigcup _{(i_1, \ldots, i_k) \in \{1,2 \ldots m\}^k}
f_{i_k} \circ \ldots \circ f_{i_1}(B_{i_1})$$

Define the \em extended atoms \em of generation $k\geq 1$ for $F_{\epsilon} $ that form $K^+$ as $f_{i_k, \epsilon} \circ \ldots \circ f_{i_1, \epsilon}(U_{i_1})$, where $(i_1, i_2, \ldots, i_k)$ is a word of length $k$ formed by symbols in $\{1, 2, \ldots, m\}$.

The diameter of each extended atom is
 smaller that $diam (B) \cdot \cdot {\lambda'} ^k$. Therefore, for sufficiently large $k\geq 1$ all
the extended atoms of generation $k$ that form $K^+$ have diameters  smaller that $\epsilon_1 /2 $.

We assert that the extended atoms of generation $k \geq 1$ are pairwise disjoint: in fact,  for two different $i \neq j $ the images  are disjoint: $f_{i, \epsilon}(U_i) \cap f_{j, \epsilon} (U_j)= \emptyset$. So the atoms of generation 1 are pairwise disjoint. Two extended atoms of generation $k$ are $f_{i_k, \epsilon} \circ \ldots \circ f_{i_1, \epsilon}(U_{i_1})$ and $f_{j_k, \epsilon} \circ \ldots \circ f_{j_1, \epsilon}(U_{j_1})$. They can intersect if and only if $(i_1, i_2, \ldots, i_k) =(j_1, j_2, \ldots, j_k)$ because each $f_{i, \epsilon}$  is an homeomorphism onto its image. So, they intersect if and only if they coincide.

By construction, $U_i \supset B_i$ and $f_{i, \epsilon}|_{B_i} = f_i$. Therefore each of the  atoms of generation $k$ for $F$, is contained in the respective extended atom of generation $k$ for $F_{\epsilon}$, that has the same finite word $(i_1, i_2, \ldots, i_k) $.

If none of the extended  atoms of generation $k$ intersects $S$, then none of the atoms of generation $k$ for $F$ intersects $S$, and the system verifies the hypothesis of Lemma \ref{lema}. So, in this case, there is nothing to prove, because the given system $F$ verifies the thesis of the Lemma \ref{lema} and thus, it is finally periodic with persistent limit cycles.

If some of
the extended atoms of generation $k$ intersects $S$, consider a new finite partition ${\cal
Q}= \{C_i\}_{1 \leq i \leq m}$ of $B$ such that the distance, defined in (\ref{definicionHausdorfDistance}),
between ${\cal Q}$ and the given partition ${\cal P }$ of $F$, is smaller
than $\epsilon_1 >0$: $$ \dist ({\cal P}, \; {\cal Q} ) < \epsilon_1 \leq \epsilon$$

We shall besides choose the new partition ${\cal Q}$ such that the new separation line  $S_{\cal Q} = \bigcup _{i
\neq j}
 (C_i \cap C_j) $ does not intersect the extended atoms of generation $k$ of $K^+$.  This last condition is possible because the diameters of the generalized atoms are all smaller than $\epsilon_1 /2$, they are compact pairwise disjoint sets,  and the distance between the two partitions ${\cal P}$ and ${\cal Q}$ (which is smaller than $\epsilon_1>0$)  can be chosen larger than $\epsilon_1 /2$, defined in (\ref{definicionHausdorfDistance}) as the maximum  Hausdorff distance between their respective pieces. (We note that the old, and principally the new, separation lines $S_{{\cal P}}$ and $S_{{\cal Q}}$, are not necessarily  $C^1$ nor even Lipschitz manifolds in the space $B$, and even if they are, they do not need to be $C^1$ or Lipschitz near one from the other, to be near with the Hausdorff distance).

  The
first condition $\dist ({\cal P}, {\cal Q}) < \epsilon_1$, joined with the assumption $\dist (U_i^c, B_i) \geq \epsilon_1 $, where $B_i$ is the $i-$th piece of the partition ${\cal P}$, implies that the respective piece $C_i$ of the partition ${\cal Q}$ verifies $C_i \subset U_i$. Therefore the
extension $f_{i, \epsilon}: U_i \mapsto B$ in $ F_{\epsilon}$   can be restricted to $C_i$.

 Define $G  = \{g_i: C_i \mapsto B\}_{1 \leq  m}$ where $g_i =
 f_{i, \epsilon}|_{C_i}$. By construction $G$ and $F$ coincide in $C_i \cap B_i$, the distance between the respective partitions ${\cal P}$ and ${\cal Q}$ is smaller than $\epsilon_1 \leq \epsilon$, and the difference of their respective contraction rates $\lambda '$ and $\lambda$ is also smaller than $\epsilon$. So $G$ is a $\epsilon  $-perturbation
 of the given $F$, according to the Definition  \ref{topologia}.  It is enough to prove now that $G$ is finally periodic with persistent limit cycles.

 Consider the limit set $K_G$ of $G$ as follows:
 $$K_{G} = \bigcap_{k \geq 1} \; \; \;  \bigcup _{(i_1, \ldots, i_k) \in \{1,2 \ldots m\}^k}
g_{i_k} \circ \ldots \circ g_{i_1}(C_{i_1})$$

As $G$ is a restriction of $F_{\epsilon } $ to the sets $C_i \subset U_i$, we
have  that $K_G \subset K^+$, and  in particular for all $k \geq 1$ the atoms of generation $k$ for $G$, i.e. $g_{i_k} \circ \ldots \circ g_{i_1}(C_{i_1})$, are contained in the extended atoms of generation $k$ for $F_{\epsilon}$.

By construction the separation line $S_{G} = S_{\cal Q}$ among the continuity pieces $C_i$ of $G$ is disjoint with the extended atoms of generation $k$ of $F_{\epsilon}$. Therefore, it is also disjoint with the atoms  of generation $k$ of $G$. Then  $G^k(B) \bigcap S_G =
\emptyset$ and, applying lemma \ref{lema}, $G$ is finally
periodic with persistent limit cycles. $\Box$

\section{Open mathematical questions.}

It is possible (but not immediate) to construct, in a  compact ball $B$ of any dimension $n-1 \geq 2$, piecewise continuous systems, uniformly locally contractive and with the separation property, as defined in Section \ref{definicionesabstractas}, that do not verify the thesis of the Theorem \ref{Teorema1}, and thus their limit set is not composed only by periodic limit cycles.

Suppose that the system had $C^r$  regularity, with $r = Lipschitz$ or with $r \geq 1$, i.e. the continuity pieces $B_i$ and the separation lines that form $S$, are bi-Lipschitz homeomorphic (or $C^r$ diffeomorphic respectively) to  $(n-1)$ or $(n-2)$ dimensional balls or manifolds, and  the homeomorphisms $f_i$ in their continuity pieces $B_i $, are bi-Lipschitz (or $C^r$-diffeomorphisms
respectively).

With this additional assumption of regularity of $F$, it is an open question to construct examples that do not verify the thesis of Theorem \ref{Teorema1}. In other words, assuming more regularity, it is unknown if the system has to  exhibit   a limit set  always formed only by periodic orbits.

On the other hand, it is also unknown if the existence of a finite number of limit cycles, as in Theorem \ref{Teorema1}, is generic for  Lipschitz or $C^1$   regular, locally contractive and piecewise continuous systems with the separation property.

\section{Conclusions}

The discontinuities of the Poincar\'{e} transformation $F$, due to spike phenomena in the neural network, play an essential role to study these systems, although it is an obstruction to apply mostly previously known results of the Theory of Dynamics Systems, which is mostly developed for continuous dynamics.

In the generic stable case, the recurrent orbits are all periodic, and all the initial states lead to limit cycles.

Due to the non-genericity of the  bifurcating case, which is a consequence of Theorem \ref{Teorema1}, those dynamic would never be seen in experiments: in fact,  arbitrarily small perturbations in the parameters of the system will lead it to periodic or quasi- periodic dynamics.  These perturbations stabilize  the system, to exhibit  a  limit set composed only by limit cycles.

Also we showed that the inter-spike interval is bounded away from zero. It means that, generically, when the system is periodic, in spite of having preferred periodic patrons of discharges, the neurons do not synchronize in phase.

On the other hand,   if the number of neurons in the system is very large, the limit cycles of the network may have a very large  period $p$, much larger than the observation time, or even than the life time of the biological system. Therefore, in spite of being asymptotically  periodic, these systems  may never show its regularity.  These two facts: extremely large periods, and irregularity inside the period,  allow us to assert  that those persistent systems with very large period $p$ shall be in fact non-predictible for the experimenter, and will be perceived   as virtually chaotic.

\end{document}